\documentclass[smallextended]{svjour3}       % onecolumn (second format)
\smartqed  % flush right qed marks, e.g. at end of proof

\makeatletter
\def\cl@chapter{\@elt {theorem}}
\makeatother

\usepackage{amssymb}
\usepackage{amsmath}
\usepackage{cite}
\usepackage[colorlinks,linkcolor=blue,citecolor=green]{hyperref}
\usepackage{cleveref}
\usepackage{microtype}

\crefname{equation}{}{}

\newtheorem{notation}[definition]{Notation}

\numberwithin{equation}{section}
\numberwithin{theorem}{section}
\numberwithin{lemma}{section}
\numberwithin{definition}{section}
\numberwithin{corollary}{section}
\numberwithin{remark}{section}
\numberwithin{proposition}{section}
\journalname{Archive for Rational Mechanics and Analysis}
\begin{document}

\title{Regularity for fully nonlinear elliptic equations with oblique boundary conditions}
\author{Dongsheng Li \and Kai Zhang$^{\ast}$
}
\institute{
          Dongsheng Li \at
              School of Mathematics and Statistics, Xi'an Jiaotong University, Xi'an 710049, China \\
              \email{lidsh@mail.xjtu.edu.cn}
\and
          Kai Zhang (*Corresponding author) \at
              School of Mathematics and Statistics, Xi'an Jiaotong University, Xi'an 710049, China; \\
              Current address: Department of Applied Mathematics, Northwestern Polytechnical University, Xi'an, Shaanxi, 710129, PR China\\
              Tel.: +86-15091636261\\
              \email{zkzkzk@stu.xjtu.edu.cn; Current email: zhang$\_$kai@nwpu.edu.cn}
}

\date{Received: date / Accepted: date}

\maketitle

\begin{abstract}
In this paper, we obtain a series of regularity results for viscosity solutions of fully nonlinear elliptic equations with oblique derivative boundary conditions. In particular, we derive the pointwise  $C^{\alpha}$, $C^{1,\alpha}$ and $C^{2,\alpha}$ regularity. As byproducts, we also prove the A-B-P maximum principle, Harnack inequality, uniqueness and solvability of the equations.
\end{abstract}

\section{Introduction}\label{S1}
We study the regularity of viscosity solutions to the following fully nonlinear elliptic equation with oblique boundary condition:
\begin{equation}\label{e1.1}
\left\{\begin{aligned}
 &F(D^2u)=f &&\mbox{in}~~\Omega; \\
 &\beta\cdot Du+\gamma u=g &&\mbox{on}~~\Gamma\subset\partial\Omega,
\end{aligned}
\right.
\end{equation}
where $\Omega \subset R^n $ is a bounded domain; $\Gamma \in C^1$ and is relatively open to $\partial\Omega$; $f$ is a function defined in $\Omega $; and $\beta $ (vector valued), $\gamma $ and $g$ are functions defined on $\Gamma$. Here ``Oblique'' means that $|\beta \cdot \vec{n}|\geq\delta_0>0$ on $\Gamma$, where $\vec{n}$ denotes the inner normal of $\Omega $. We call $\Gamma$ the oblique boundary. Since the sign of $\gamma$ is not required, without loss of generality, throughout this paper, we always assume that
\begin{equation*}
  \beta\cdot \vec{n}\geq\delta_0 ~\mbox{on}~\Gamma \mbox{ and } \|\beta\|_{L^{\infty }(\Gamma) }\leq 1
\end{equation*}
for some positive constant $\delta_0$. In addition, we always assume that $F$ is uniformly elliptic, i.e., there exist positive constants $\lambda $ and $\Lambda $ such that
\begin{equation*}
  \lambda\|N\|\leq F(M+N)-F(M)\leq \Lambda \|N\|, ~~~\forall M,N\in S^n,~N\geq 0,
\end{equation*}
where $S^n$ denotes the set of $n\times n$ symmetric matrices; $N\geq 0$ means the nonnegativeness and $\|N\|$ is the spectral radius of $N$. Since we deal with the viscosity solutions of\cref{e1.1}, it is convenient for us to assume that $f$, $\beta$, $\gamma$ and $g$ are continuous throughout this paper.

Few regularity but some a priori estimates and existence results are known for\cref{e1.1}, where $\Gamma=\partial\Omega$ and $\gamma\leq0$ are needed. For linear equations, the existence of $C^{2,\alpha}$ solutions are obtained by Green's representation, the method of freezing the coefficients and method of continuity (cf. \cite[Chapter 6]{G-T}). The work of Lieberman \cite{Li2} covered above results and only required the $C^{1,\alpha}$ smoothness on $\partial \Omega$. Safonov \cite{Sa1, Sa2} extended this results to the nonlinear Bellman equations, that is, he proved the existence of $C^{2,\alpha}$ solutions as $\partial\Omega\in C^{1,\alpha}$ for\cref{e1.1}.
 In 1982, Lieberman \cite{Li1} obtained the existence of solutions for quasilinear equations based on a fixed point theorem and the solvability for linear equations. By a priori estimates and the method of continuity, Lieberman and Trudinger \cite{L-T} proved the existence and uniqueness of $C^{2,\alpha}$ solutions for fully nonlinear elliptic equations with nonlinear oblique boundary conditions as $\partial\Omega\in C^4$.

Since 1980s, the notion of viscosity solutions has been applied widely in the study of non-divergent equations especially of fully nonlinear elliptic equations; and some important interior
regularity results and global results with Dirichlet boundary conditions have been obtained (see \cite{C-C} and \cite{C-I-L} and references therein). As for applications of viscosity solutions to oblique derivative problems, Ishii \cite{Is1} obtained the existence and uniqueness of viscosity solutions for fully nonlinear elliptic equations. Milakis and Silvestre \cite{M-S} proved the $C^{1,\alpha}$ and $C^{2,\alpha}$ regularity results for fully nonlinear equations with Neumann boundary data on flat boundaries. For more a priori estimates of the oblique derivative problems, we refer the reader to the book \cite{Li7}.

In this paper, we derive a series of regularity results for viscosity solutions of\cref{e1.1}. In particular, the $C^{\alpha}$, $C^{1,\alpha}$ and $C^{2,\alpha}$ boundary regularity are deduced:

\begin{theorem}[$C^{\alpha}$ regularity]\label{th5.2}
Let $u$ satisfy
  \begin{equation*}
    \left\{
    \begin{aligned}
     & u\in S(\lambda,\Lambda,f) &&\mbox{in}~~ \Omega ; \\
     & \beta \cdot Du +\gamma u= g &&\mbox{on}~~\Gamma.
    \end{aligned}
    \right.
  \end{equation*}

Then for any $\Omega'\subset\subset \Omega\cup\Gamma$, $u\in C^{\alpha}(\bar{\Omega'})$ and
\begin{equation*}
  \|u\|_{C^{\alpha}(\bar{\Omega'})}\leq C\left(\|u\|_{L^{\infty }(\Omega)}+\|f\|_{L^n(\Omega)}+\|g\|_{L^{\infty }(\Gamma )}\right),
\end{equation*}
where $0<\alpha <1$ depends only on $n$, $\lambda$, $\Lambda$ and $\delta_0$, and $C$ depends also on $\|\gamma\|_{L^{\infty }(\Gamma )}$, $\Omega '$ and $\Omega $.
\end{theorem}

\begin{theorem}[$C^{1,\alpha}$ regularity]\label{co4.1}
Let $u$ be a viscosity solution of\cref{e1.1} and $0<\alpha <\bar{\alpha}$ where $0<\bar{\alpha}<1$ is a constant depending only on $n$, $\lambda$, $\Lambda$ and $\delta_0$. Suppose that $\beta, \gamma, g\in C^{\alpha }(\bar{\Gamma} )$ and $f$ satisfies
\begin{equation*}
\left (\frac{1}{|B_r(x_0)\cap\Omega|}\int_{B_r(x_0)\cap\Omega}|f|^n\right )^{\frac{1}{n}}\leq C_f r^{\alpha - 1}~~\forall x_0\in \Omega,~~\forall r>0.
\end{equation*}

Then for any $\Omega' \subset \subset \Omega\cup \Gamma$, $u\in C^{1,\alpha}(\bar{\Omega '})$ and
\begin{equation}\label{e4.42}
   \|u\|_{C^{1,\alpha }(\bar{\Omega '})}\leq C\left(\|u\|_{L^{\infty }(\Omega )}+C_f+\|g\|_{C^{\alpha }(\bar{\Gamma}  )}+|F(0)|\right),
\end{equation}
where $C$ depends only on $n$, $\lambda$, $\Lambda$, $\delta_0$, $\alpha $, $\|\beta \|_{C^{\alpha }(\bar{\Gamma} )}$, $\|\gamma  \|_{C^{\alpha }(\bar{\Gamma })}$, $\Omega '$ and $\Omega$.
\end{theorem}

\begin{theorem}[$C^{2,\alpha}$ regularity]\label{th4.6}
Let $F$ be convex, $u$ be a viscosity solution of\cref{e1.1} and $0<\alpha <\tilde{\alpha}$ where $0<\tilde{\alpha}<1$ is a constant depending only on $n$, $\lambda$, $\Lambda$ and $\delta_0$. Suppose that $\Gamma\in C^{1,\alpha} $, $\beta,\gamma,g\in C^{1,\alpha }(\bar{\Gamma}  )$ and $f\in C^{\alpha }(\bar{\Omega}  )$.

Then for any $\Omega' \subset \subset \Omega\cup \Gamma$, $u\in C^{2,\alpha}(\bar{\Omega '})$ and
\begin{equation*}
   \|u\|_{C^{2,\alpha }(\bar{\Omega '})}\leq C\left(\|u\|_{L^{\infty }(\Omega  )}+\|f\|_{C^{\alpha }( \bar{\Omega} )}+\|g\|_{C^{1, \alpha }(\bar{\Gamma}  )}+|F(0)|\right),
\end{equation*}
where $C$ depends only on $n$, $\lambda$, $\Lambda$, $\delta_0$, $\alpha $, $\|\beta \|_{C^{1,\alpha }(\bar{\Gamma}  )}$, $\|\gamma  \|_{C^{1,\alpha }(\bar{\Gamma}  )}$, $\Omega '$and $\Omega$.
\end{theorem}

The following is the outline of this paper. First, the Alexandrov-Bakel'man-Pucci type maximum principle and boundary Harnack type inequality are presented as basic tools. Then $C^{\alpha}$ regularity follows clearly. The $C^{1,\alpha}$ and $C^{2,\alpha}$ regularity for viscosity solutions of\cref{e1.1} are obtained by approximating the original problem by model problems at each scales whose regularity can be dealt with and which should be designed properly. We point out here that we do not need flatten the boundary, which is different from the references mentioned above. This paper also contains uniqueness and existence results of\cref{e1.1}, which will be used to solve our model problems.

The paper is organized in the following way. We introduce the Alexandrov-Bakel'man-Pucci (A-B-P for short) maximum principle and the boundary Harnack inequality in \Cref{S2}, which are the basic tools to attack the regularity for\cref{e1.1}. Based the boundary Harnack inequality, the pointwise $C^{\alpha}$ regularity follows. The proof for the A-B-P maximum principle has been motivated by \cite{M-S}, where the authors deal with the Neumann problems. Combining the A-B-P maximum principle, the interior Harnack inequality and the barrier technique, we derive the boundary Harnack inequality. The barrier is adopted from \cite{L-T}.

A Jensen's type uniqueness of viscosity solutions is presented in \Cref{S3}. That is, the subsolution minus the supersolution is also a ``subsolution''. This leads to the uniqueness of viscosity solutions by combining with the A-B-P maximum principle, where
a mixed boundary value problem is considered. In addition, we also prove the existence result by Perron's method, which will be used to construct our model problem to approximate the viscosity solution of\cref{e1.1} later.

We present our model problem and prove the $C^{1,\alpha}$ and $C^{2,\alpha}$ regularity of its solution in \Cref{S4}, which will be used to approximate the viscosity solution of\cref{e1.1} at different scales. Roughly speaking, the model problem is that $f,\gamma,g\equiv 0$, $\beta$ is a constant vector and $\Omega$ is a spherical cap with $\Gamma$ is the flat part. The $C^{1,\alpha}$ regularity is a consequence of the Jensen's type uniqueness obtained in \Cref{S3}. To show the $C^{2,\alpha}$ regularity under the assumption that $F$ is convex, we first prove that the restriction of the viscosity solution on the flat boundary is also a viscosity solution in dimension $n-1$, and then by the interior $C^{2,\alpha}$ estimates for convex operators, we have the $C^{2,\alpha}$ regularity of the restriction solution. Finally, from the boundary $C^{2,\alpha}$ regularity for the Dirichlet problem, we deduce the conclusion.

Based on the previous results, we show the pointwise $C^{1,\alpha}$ and
$C^{2,\alpha}$ regularity for\cref{e1.1} in \Cref{S5} and  \Cref{S6} respectively, where higher regularities for\cref{e1.1} is also presented in \Cref{S6}.
Our proof is influenced by the spirits in \cite{C-C}.

Before the end of this section, we introduce the following notations which will be frequently used in this paper.
\begin{notation}\label{no1.1}
\begin{enumerate}
\item $\{e_i\}^{n}_{i=1}$: the standard basis of $R^n$, i.e., $e_i=(0,...0,\underset{i^{th}}{1},0,...0)$.
\item Given $x=(x_1,x_2,...,x_n)\in R^{n}$, we may rewrite $x=(x',x_{n})$ where $x'=(x_1,x_2,...,x_{n-1})$ .
\item $S^n$: the set of $n\times n$ symmetric matrices and $\|A\|:=$ the spectral radius of $A$ for any $A\in S^n$.
\item $R^n_+:=\{x\in R^n\big|x_n>0\}$.
\item Given $\beta \in R^n $ with $\beta_n\neq 0$, denote by $x_{\beta }'=(x_{\beta,1 },x_{\beta,2 },...,x_{\beta,n-1 })$ the projection of $x$ to the hyperplane $\{x_n=0\}$ along the direction $\beta$, i.e., $x_{\beta, i }=x_i-(\beta _i/\beta _n) x_n$ for $i=1,2,...,n-1$. Clearly, $|x-(x_{\beta }',0)|=(|\beta |/|\beta _n|) |x_n|$.
\item $B_r(x_0):=\{x\in R^{n}\big| |x-x_0|<r\}$, $B_r:=B_r(0)$ and $B_r^+(x_0):=B_r(x_0)\cap R^n_+$.
\item $T_r(x_0)\ :=\{(x',0)\in R^{n}\big| |x'-x_0'|<r\}$ \mbox{ and } $T_r:=T_r(0)$.
\item $B_{r,h}:=B_R(-(R-h)e_n)$, where $R$ satisfies $(R-h)^2+r^2=R^2$.
\item $B^+_{r,h}$: the spherical cap with the base radius $r$ and height $h$ ($h\leq r$), i.e., $B^+_{r,h}:=B_{r,h}\cap R^n_+$.
\item $A^c:$ the complement of $A$, $\overset{\circ}{A}:$ the interior of $A$ and $\bar A: $ the closure of $A$, $\forall A\subset R^n$.
\item $\mathrm{dist}(A,B):=\inf\{|x-y|\big|x\in A, y\in B\}$, $~\forall~A,B\subset R^n$.
\item $a^{+}:=\max\{0,a\}$ \mbox{ and } $a^{-}:=-\min\{0,a\}$.
\item $\vec{n}(x_0):$ the inner normal of $\Omega $ at $x_0\in\partial\Omega$.
\item $I:$ the unit matrix.
\item Given a function $\varphi $. We may use $\varphi _i$ or $D_i \varphi $ to denote $\partial \varphi/\partial x _{i}$. Similarly, $\varphi _{ij}$ and $D_{ij}\varphi $ denote $\partial ^{2}\varphi/\partial x_{i}\partial x_{j}$.
\item $D \varphi:=(\varphi_1 ,...,\varphi_{n} )$, $D_{\_}\varphi:=( \varphi_1 ,..., \varphi_{n-1} )$, $D^2 \varphi :=\left(\varphi _{ij}\right)_{n\times n}$ and $D^2_{\_}\varphi:=\left(\varphi _{ij}\right)_{(n-1)\times (n-1)}$.
\item For $k\geq 0$, $\varphi $ is called $C^{k,\alpha}$ at $x_0$ if there exists a polynomial $P$ of degree $k$ such that
$|\varphi (x) - P (x)|\leq K |x-x_0|^{k+\alpha}$ for any $x$. Then, $D^{\zeta}\varphi (x_0):=D^{\zeta } P(x_0)$ where $\zeta $ is a multi index, $[\varphi ]_{C^{k,\alpha}(x_0)}:=K$ and $\|\varphi \|_{C^{k,\alpha}(x_0)}:=K+\sum_{|\zeta|\leq k}|D^{\zeta }\varphi (x_0)|  $.
\item For viscosity solutions, we use the notations $\bar{S}(\lambda ,\Lambda ,f)$, $\underline{S}(\lambda ,\Lambda ,f)$, $S(\lambda ,\Lambda ,f)$, $M^{+}(M,\lambda ,\Lambda )$, $M^{-}(M,\lambda ,\Lambda )$ etc. as in \cite{C-C}.
\end{enumerate}
\end{notation}

\section{A-B-P maximum principle, Harnack inequality and \texorpdfstring{$C^{\alpha}$}{C, alpha} regularity}\label{S2}
In this section, we introduce some notations and present some elementary results concerning the viscosity solutions of oblique derivative problems. We say that $u$ touches $v$ by above at $x_0$ in $\Omega $ if $u\geq v$ in $\Omega$ and $u(x_0)=v(x_0)$. Similarly, we have the definition for touching by below. Now, we give the definition of viscosity solutions.
\begin{definition}\label{de1.1}
Let $u$ be upper (lower) semicontinuous in $\Omega \cup \Gamma  $. We say that $u$ is a viscosity subsolution (supersolution) of\cref{e1.1} if for any $\varphi\in C^2(\Omega \cup \Gamma)$ touching $u$ by above (below) at $x_0$ in $\Omega \cup \Gamma $, we have that
\begin{equation*}
\begin{aligned}
F(D^2\varphi (x_0))\geq (\leq) f(x_0) ~~\mbox{if}~~ x_0\in\Omega
\end{aligned}
\end{equation*}
and
\begin{equation}\label{e6.2}
\begin{aligned}
\beta (x_0)\cdot D\varphi (x_0)+\gamma(x_0) \varphi (x_0) \geq (\leq) g(x_0) ~~\mbox{if}~~ x_0\in\Gamma.
\end{aligned}
\end{equation}

If $u\in C(\Omega \cup \Gamma)$ is both subsolution and supersolution, we call it a viscosity solution.
\end{definition}
\begin{remark}\label{r2.2}
(i) We may write $u\in USC(\Omega \cup \Gamma)$ ($u\in LSC(\Omega \cup \Gamma)$) for short if $u$ is upper (lower) semicontinuous in $\Omega \cup \Gamma $.

(ii) Touching in $\Omega \cup \Gamma$ can be replaced by touching in a neighborhood of $x_0$ (see \cite[Proposition 2.4]{C-C}). In this case, $\varphi \in C^2(\Omega \cup \Gamma)$ in the definition can be replaced by that $\varphi$ is a paraboloid (a polynomial of degree 2).

(iii) The notion of viscosity solutions for oblique derivative problems was introduced first by P.-L. Lions \cite{Lio}.
\end{remark}

Now we present a closedness result concerning the viscosity solutions.
\begin{proposition}\label{le2.3}
Suppose that $\Gamma\in C^2$. Let $\{u_m\}$ satisfy
  \begin{equation*}
    \left\{
    \begin{aligned}
     &u_m\in S(\lambda,\Lambda,f) &&\mbox{in}~~ \Omega; \\
     &\beta \cdot Du_m +\gamma u_m = g &&\mbox{on}~~\Gamma.
    \end{aligned}
    \right.
  \end{equation*}
Suppose that $u_m\rightarrow u$ uniformly in any compact subset of $\Omega \cup \Gamma$. Then
  \begin{equation*}
    \left\{
    \begin{aligned}
     &u\in S(\lambda,\Lambda,f) &&\mbox{in}~~ \Omega; \\
     &\beta \cdot Du +\gamma u = g &&\mbox{on}~~\Gamma.
    \end{aligned}
    \right.
  \end{equation*}
\end{proposition}
\proof  $u\in S(\lambda ,\Lambda,f )$ in $\Omega $ is classical (see Proposition 2.9 \cite{C-C}). We only prove that $\beta\cdot Du+\gamma u\geq g$ on $\Gamma$ and the proof for the other direction is similar. Let $P$ be a paraboloid touching $u$ by above in a neighborhood of $x_0\in \Gamma$. Without loss of generality, we assume that $x_0=0$ and $\vec{n}(0)=e_n$. Let $\varphi $ denote the representation function of $\Gamma$ near $0$ with $\varphi (0)=0$ and $D\varphi (0)=0$. We need to prove that $\beta (0)\cdot DP(0)+\gamma(0)P(0)\geq g(0)$. Suppose not. Then
\begin{equation}\label{e1zk.36}
  P_n(0)<\frac{1}{\beta _{n}(0)}\left(g(0)-\gamma(0)P(0)-\beta'(0)\cdot D_{\_}P(0)\right):=A.
\end{equation}
By the continuity of $P_n$, there exist $\tau,\eta>0$ such that
\begin{equation*}
  P_n(x)<A-\eta~~\forall x\in B_{\tau}.
\end{equation*}
By Taylor's formula, for any $x\in\bar{\Omega}\cap B_\tau$,
\begin{equation*}
\begin{aligned}
    P (x)=&P (x',\varphi(x'))+P _{n}(x',\varphi(x'))(x_n-\varphi(x'))+\frac{1}{2}P _{nn}(x_n-\varphi(x'))^2\\
\leq &P (x',\varphi(x'))+(A-\eta)(x_n-\varphi(x'))+\frac{1}{2}P_{nn}(x_n-\varphi(x'))^2.
\end{aligned}
\end{equation*}
By the boundedness of $P_{nn} $, for any $N>0$ there exists $\tau '> 0$ such that
\begin{equation*}
P (x)\leq P (x',\varphi(x'))+(A-\eta/2)(x_n-\varphi(x'))-N(x_n-\varphi(x'))^2 ~~\forall x\in \bar{\Omega}\cap B_{\tau'}.
\end{equation*}
The constant $N$ is to be chosen later. Here, we choose $\tau '$ small enough such that $\Omega \cap B_{\tau'}\subset \subset \Omega \cup \Gamma$ and
\begin{equation}\label{e4-14.1}
  \beta_n\geq \delta_0/2 ~\mbox{on}~ \Gamma\cap B_{\tau'}.
\end{equation}
For $ \varepsilon > 0$, let $\psi_{\varepsilon } (x)=P (x',\varphi(x'))+(A-\eta/2)(x_n-\varphi(x'))-N(x_n-\varphi(x'))^2+\varepsilon |x|^2$. Then
\begin{equation*}
  \psi_{\varepsilon } \geq P \geq u ~~\mbox{in}~~\Omega \cap B_{\tau'}\mbox{ and }\psi_{\varepsilon } (0)=P (0)=u(0).
\end{equation*}

Since $u_m\rightarrow u$ uniformly in $\Omega \cup \Gamma$, there exists $m_0$ large enough such that
\begin{equation*}
  \|u_{m_0}-u\|_{L^{\infty }(\Omega \cap B_{\tau'})}\leq \frac{1}{2}(\varepsilon \tau ')^{2}.
\end{equation*}
Then
\begin{equation*}
  |\psi_{\varepsilon } (0)-u_{m_0} (0)|\leq \frac{1}{2}(\varepsilon \tau ')^{2}~~\mbox{and}~~\psi_{\varepsilon } -u_{m_0}> \frac{1}{2}(\varepsilon \tau ')^2~~\mbox{on}~~\partial (\Omega \cap B_{\tau'})\backslash \Gamma.
\end{equation*}
Hence, $\psi_{\varepsilon } +c_0$ touches $u_{m_0}$ by above at some $\tilde{x}\in \Omega \cap B_{\tau'}$ for a proper $c_0$ with $|c_0|\leq (\varepsilon \tau ')^{2}/2.$

By a calculation, the second derivatives of $\psi_{\varepsilon}$ are
\begin{equation*}
\psi_{\varepsilon ,ij}=\left\{
  \begin{aligned}
  &(-A+\eta/2-2N\varphi+2Nx_{n})\varphi_{ij}-2N\varphi _{i}\varphi _{j} \\
  &+P_{ij}+P_{in}\varphi _{j}+P_{jn}\varphi _{i}+P_n\varphi_{ij}+2\varepsilon , && i,j<n,\\
  &2N\varphi _{i}, && i<n,j=n,\\
  &-2N+2\varepsilon, && i=j=n.
  \end{aligned}\right.
\end{equation*}
Note that $\varphi (0)=0$ and $D\varphi (0)=0$. Then by choosing $N$ large enough and $\tau'$ small enough, we have
\begin{equation*}
  M^+\left(D^2\psi_{\varepsilon}\right)<-\|f\|_{L^{\infty }(\Omega )}~~\mbox{in}~~\Omega\cap B_{\tau'}.
\end{equation*}
Combining with $M^{+}\left(D^2u_{m_0}\right)\geq f$, we know that $\tilde{x}\notin {\Omega}$, i.e., $\tilde{x}\in \Gamma$. By the definition of viscosity solutions, we have
\begin{equation*}
\beta (\tilde{x})\cdot D\psi_{\varepsilon } (\tilde{x})+\gamma(\tilde{x}) \psi_{\varepsilon } (\tilde{x}) \geq g(\tilde{x}).
\end{equation*}
By a calculation (note that $\tilde{x}_n=\varphi (\tilde{x}')$),
\begin{equation*}
\begin{aligned}
  D\psi_{\varepsilon }=(A-\eta/2)(-D\varphi,1)- 2\varepsilon \tilde{x}+(D_{\_}P+P_nD\varphi,0)~~\mbox{at}~~ \tilde{x}.
\end{aligned}
\end{equation*}
Hence, at $\tilde{x}$
\begin{equation*}
g\leq (\eta/2-A)\beta' D\varphi+\beta _n(A-\eta/2)-2\varepsilon \beta \tilde{x}+\beta ' (D_{\_}P+P_nD\varphi)+\gamma u.
\end{equation*}
Note that $D\varphi(0)=0$ and $|\tilde{x}|\leq \sqrt{\varepsilon}\tau '$. Taking $\varepsilon $ small enough, we have
\begin{equation*}
g\leq \beta _n(A-\eta/4)+\beta ' D_{\_}P+\gamma u~~\mbox{at}~~ \tilde{x}.
\end{equation*}
By the definition of $A$,
\begin{equation*}
\begin{aligned}
  g(\tilde{x})= &\frac{\beta _n (\tilde{x})}{\beta _n (0)}\left(g(0)-\gamma(0)P (0)-\beta'(0)\cdot D_{\_} P(0)\right)-\beta _n (\tilde{x})\eta/4\\
  &+\beta' (\tilde{x})\cdot D_{\_}P (\tilde{x})++\gamma(\tilde{x}) u (\tilde{x}).
\end{aligned}
\end{equation*}
By the continuity of each ingredient in above equation and taking $\varepsilon $ small enough, we have
\begin{equation*}
  -\beta _n (\tilde{x})\eta/8\geq 0
\end{equation*}
which is impossible by recalling\cref{e4-14.1}.\qed~\\

Next, we intend to present an Alexandrov-Bakel'man-Pucci type maximum principle for oblique derivative problems following the idea of \cite{M-S}.
\begin{lemma}\label{th2.4}
Let $\Omega\subset B_1$ and $u$ satisfy
  \begin{equation*}
    \left\{
    \begin{aligned}
     &u\in \bar{S}(\lambda,\Lambda,f) &&\mbox{in}~~ \Omega; \\
     &\beta \cdot Du \leq g &&\mbox{on}~~\Gamma.
    \end{aligned}
    \right.
  \end{equation*}
Suppose that there exists a direction $\xi\in\partial B_1$ such that $\beta\cdot \xi \geq \delta_1$ on $\Gamma$. Then
\begin{equation}\label{e1.29}
  \sup_{\Omega}u^{-}\leq \sup_{\partial \Omega \backslash \Gamma}u^-+C \max_{\Gamma}g^++C\|f^+\|_{L^n(\{u=\Gamma_u\})},
\end{equation}
where $\Gamma_u$ is the convex envelop of $u$ in $\Omega$ and $C$ depends only on $n$, $\lambda$, $\Lambda$ and $\delta_1$.
\end{lemma}
\proof  Without loss of generality, let $\xi$ be $e_n$. We assume that $u \geq 0$ on $\partial \Omega \backslash \Gamma$. Otherwise, we may consider $\tilde{u}:=u-\underset{\partial \Omega \backslash \Gamma}{\inf}u$. Let $M:=\underset{\Omega }{\sup}~u^{-}$ and
\begin{equation}\label{e1.30}
  \textbf{A}:=\left\{A\in R^n \big|A_n\delta_1>2\max_{\Gamma}g,|A'|\leq \frac{\delta_1 }{2}A_n \mbox{ and } |A|\leq \frac{M}{4}\right\}.
\end{equation}

For any $A\in \textbf{A}$, $A\cdot x+c_0$ touches $u$ by below at some $x_0\in \bar{\Omega}$ for a proper $c_0\in R$. Since $|A|\leq M/4$ and $\Omega\subset B_1$, $x_0\notin \partial \Omega \backslash \Gamma$. If $x_0\in \Gamma$, then
\begin{equation*}
  g(x_0)\geq A_n\beta_n(x_0)+A'\cdot \beta'(x_0)\geq A_n\delta_1-\frac{\delta_1 }{2}A_n= \frac{\delta_1 }{2}A_n>\max_{\Gamma}g,
\end{equation*}
which is a contradiction. Hence $x_0\in \Omega$, i.e., $A\in \nabla \Gamma_u(\Omega) $.

If $\max_{\Gamma}g>\delta_1 M/16$, then
\begin{equation}\label{e1.31}
\sup_{\Omega}u^{-}=M \leq \frac{16}{\delta_1}\max_{\Gamma}g.
\end{equation}
Otherwise,
\begin{equation*}
\left\{A\in R^n \big|A_n>\frac{M}{8},|A'|\leq \frac{\delta_1 }{2}A_n,|A|\leq \frac{M}{4}\right\}\subset \textbf{A}\subset \nabla \Gamma_u(\Omega).
\end{equation*}
From the proof of \cite[Theorem 3.2]{C-C}, we have that $\Gamma_u\in C^{1,1}(\bar{\Omega})$ and
\begin{equation*}
  \|f^+\|^n_{L^n(\{u=\Gamma_u\})}\geq C|\nabla \Gamma_u(\Omega)|,
\end{equation*}
where $C$ depends only on $n$, $\lambda$ and $\Lambda$. Hence,
\begin{equation}\label{e1.32}
  \|f^+\|^n_{L^n(\{u=\Gamma_u\})}\geq CM^n.
\end{equation}

Combining\cref{e1.31} and\cref{e1.32}, we conclude that
\begin{equation*}
  \sup_{\Omega}u^{-}=M\leq C \max_{\Gamma}g^++C\|f^+\|_{L^n(\{u=\Gamma_u\})}.
\end{equation*}
\qed ~\\

Now, the full version of A-B-P maximum principle follows easily:
\begin{theorem}[\textbf{A-B-P}]\label{th_ABP.1}
Let $\Omega\subset B_1$ and $u$ satisfy
\begin{equation*}
    \left\{
    \begin{aligned}
     &u\in \bar{S}(\lambda,\Lambda,f) &&\mbox{in}~~ \Omega; \\
     &\beta \cdot Du +\gamma u \leq g &&\mbox{on}~~\Gamma.
    \end{aligned}
    \right.
\end{equation*}
Suppose that $\gamma\leq 0$ on $\Gamma$ and there exists $\xi\in\partial B_1$ such that $\beta\cdot \xi \geq \delta_1$ on $\Gamma$. Then
\begin{equation}\label{e_ABP1.29}
  \sup_{\Omega}u^{-}\leq \sup_{\partial \Omega \backslash \Gamma}u^-+C \max_{\Gamma}g^++C\|f^+\|_{L^n(\{u=\Gamma_u\})},
\end{equation}
where $C$ depends only on $n$, $\lambda$, $\Lambda$ and $\delta_1$.
\end{theorem}
\proof  Let $v=\min\{u,0\}$. Then (note that $\gamma\leq 0$ and $v\leq 0$)
\begin{equation*}
    \left\{
    \begin{aligned}
     &v\in \bar{S}(\lambda,\Lambda,f) &&\mbox{in}~~ \Omega; \\
     &\beta \cdot Dv \leq g^+  &&\mbox{on}~~\Gamma.
    \end{aligned}
    \right.
\end{equation*}
Hence, by \Cref{th2.4}, we have
\begin{equation*}
\begin{aligned}
  \sup_{\Omega}u^{-}=\sup_{\Omega}v^{-}&\leq \sup_{\partial \Omega \backslash \Gamma}v^-+C \max_{\Gamma}g^++C\|f^+\|_{L^n(\{v=\Gamma_v\})}\\
  &= \sup_{\partial \Omega \backslash \Gamma}u^-+C \max_{\Gamma}g^++C\|f^+\|_{L^n(\{u=\Gamma_u\})}.
\end{aligned}
\end{equation*}\qed~\\

\begin{remark}\label{re.uni}
The hypothesises $\beta \cdot \vec{n}\geq \delta_0$ and $\beta\cdot \xi\geq \delta_1$ on $\Gamma$ imply that the A-B-P maximum principle holds when $\Gamma$ is a ``small'' portion of $\partial\Omega$ (see also \cite{WLH}).
\end{remark}

Through above A-B-P maximum principle, we obtain the boundary Harnack type inequality (see also \cite{L-T}).
\begin{theorem}[\textbf{Boundary Harnack inequality}]\label{t1.2}
Let $u\geq 0$ satisfy
  \begin{equation*}
    \left\{
    \begin{aligned}
     & u\in S(\lambda,\Lambda,f) &&\mbox{in}~~ \Omega; \\
     & \beta \cdot Du +\gamma u= g &&\mbox{on}~~\Gamma.
    \end{aligned}
    \right.
  \end{equation*}
Suppose that $\gamma\leq 0$ and $\beta_n \geq \delta_1$ on $\Gamma$. In addition, assume that
  \begin{equation}\label{e1.13}
\Gamma=\{(x',x_n)\in B_{1}\big|x_n=\varphi (x')\}, \{(x',x_n)\in B_{1}\big|x_n>\varphi (x')\}\subset \Omega
\end{equation}
and $\varphi $ satisfies
\begin{equation}\label{e1.14}
\varphi(0)=0~~\mbox{and}~~ D\varphi(0)=0.
\end{equation}

Then there exist constants $0<\rho<1$ and $C$ depending only on $n$, $\lambda$, $\Lambda$, $\delta_1$ and $\|\gamma \|_{L^{\infty }(\Gamma)}$, and a constant $0<R_0<1$ depending also on the $C^1$ continuity modulus of $\Gamma$ such that for any $R\leq R_0$,
\begin{equation}\label{e1.33}
  \sup_{\tilde{G}(R)}u\leq C\left(\inf_{G(R/4)}u+R\|f\|_{L^n(\Omega)}+R\|g\|_{L^{\infty }(\Gamma)}\right),
\end{equation}
where $G(R):=\left\{x\in\Omega \big||x'|<R,x_n<\rho R\right\}$ and $\tilde{G}(R):=\{x\in\Omega\big||x'|<R$,$\rho R<x_n<3\rho R\}$.
\end{theorem}
\proof  Take $\rho$ small enough such that
\begin{equation}\label{e1.43}
  \rho <\frac{\delta_1 }{16(1+\|\gamma \|_{L^{\infty }(\Gamma) })}
\end{equation}
and
\begin{equation}\label{e1.41}
  M^+ \left(\left(
                  \begin{array}{cc}
                    2I & 0 \\
                    0 & -\frac{1}{2\rho ^{2}} \\
                  \end{array}
                \right)\right)<-1.
\end{equation}
Since $\varphi (0)=0$ and $D\varphi (0)=0$, we take $R_0$ small enough such that
\begin{equation}\label{e1.42}
  |\varphi (x')|\leq \frac{\rho }{16}|x'|~~ ~~\forall |x'|\leq R_0.
\end{equation}

Then for any $R\leq R_0/2$, $\mathrm{dist}(\tilde{G}(R),\partial \Omega)\geq CR$. By the interior Harnack inequality (see \cite[Theorem 4.3]{C-C} with a proper scaling)
\begin{equation*}
   \sup_{\tilde{G}(R)}u\leq C\bigg(\inf_{\tilde{G}(R)}u+R\|f\|_{L^n(\Omega)}\bigg),
\end{equation*}
we only need to prove
\begin{equation}\label{e2.8}
\inf_{\tilde{G}(R)}u\leq C\left(\inf_{G(R/4)}u+R\|f\|_{L^n(\Omega)}+R\|g\|_{L^{\infty }(\Gamma)}\right).
\end{equation}
Let $A:=\underset{\tilde{G}(R)}{\inf}u$. Set
\begin{equation*}
  w_1(x)=2\rho R-x_n, w_2(x)=2-\left(\frac{x_n}{2\rho R}\right)^2-\frac{x_n}{2\rho R}+\frac{|x'|^2}{R^2}
\end{equation*}
and
\begin{equation*}
  w(x)=u(x)+\frac{1}{\delta_1 }\|g\|_{L^{\infty }(\Gamma)}w_1(x)+\frac{1}{4} Aw_2(x)-A.
\end{equation*}

Combining\cref{e1.43},\cref{e1.41},\cref{e1.42} with the definition of $A$, it is easy to verify that
\begin{equation*}
  \left\{
  \begin{aligned}
     & w\in \bar{S}(\lambda /n,\Lambda,f) &&\mbox{in}~~G(2R); \\
     & w\geq 0 &&\mbox{on}~~\partial G(2R)\backslash\Gamma;\\
     & \beta\cdot Dw+\gamma w\leq 0 &&\mbox{on}~~\Gamma .
  \end{aligned}
  \right.
\end{equation*}
By the A-B-P maximum principle,
\begin{equation*}
  w\geq -CR\|f\|_{L^{n}(\Omega)},
\end{equation*}
i.e.,
\begin{equation*}
  u+CR\|f\|_{L^{n}(\Omega)}+CR\|g\|_{L^{\infty }(\Gamma)}\geq A-\frac{1}{4} A w_2~~\mbox{in}~~G(2R).
\end{equation*}
Then
\begin{equation*}
  u+CR\|f\|_{L^{n}(\Omega)}+CR\|g\|_{L^{\infty }(\Gamma)}\geq \frac{ A}{4}~~\mbox{in}~~G(R/4),
\end{equation*}
which is\cref{e2.8}. \qed~\\

Based on above boundary Harnack inequality, the boundary pointwise H\"{o}lder estimate follows standardly (see \cite[Theorem 8.22 and Theorem 9.31]{G-T}):
\begin{lemma}\label{th5.1}
Let $u$ satisfy
  \begin{equation*}
    \left\{
    \begin{aligned}
     & u\in S(\lambda,\Lambda,f) &&\mbox{in}~~ \Omega ; \\
     & \beta \cdot Du= g &&\mbox{on}~~\Gamma
    \end{aligned}
    \right.
  \end{equation*}
and $x_0\in \Gamma$ such that $\mathrm{dist}(x_0,\partial\Omega\backslash \Gamma)>1$.

Then $u$ is $C^{\alpha_0}$ at $x_0$. Precisely, for any $r\leq \check{C}^{-1}$,
\begin{equation}\label{e5.100}
  \|u-u(x_0)\|_{L^{\infty }(\bar{\Omega }\cap B_r(x_0))}\leq C\left(\|u\|_{L^{\infty }(\Omega)}+\|f\|_{L^n(\Omega)}+\|g\|_{L^{\infty }(\Gamma )}\right)r^{\alpha_0},
\end{equation}
where $0<\alpha_0 <1$ and $C$ depend only on $n$, $\lambda$, $\Lambda$ and $\delta_0$, and $\check{C}$ depends also on the $C^{1}$ modulus of $\Gamma$ at $x_0$.
\end{lemma}
\begin{remark}\label{re5.1}
The condition $\mathrm{dist}(x_0,\partial\Omega\backslash \Gamma)>1$ is not an essential assumption and ``1'' can be replaced by any positive constant. Then, we obtain the scaling version of\cref{e5.100}.
\end{remark}

Then we have the following pointwise $C^{\alpha}$ estimate:
\begin{theorem}\label{th_ca.1}
Let $u$ satisfy
  \begin{equation*}
    \left\{
    \begin{aligned}
     & u\in S(\lambda,\Lambda,f) &&\mbox{in}~~ \Omega ; \\
     & \beta \cdot Du +\gamma u= g &&\mbox{on}~~\Gamma
    \end{aligned}
    \right.
  \end{equation*}
and $x_0\in \Gamma$ such that $\mathrm{dist}(x_0,\partial\Omega\backslash \Gamma)>1$.

Then $u$ is $C^{\alpha_0}$ at $x_0$. Precisely, for any $r\leq \check{C}^{-1}$,
\begin{equation}\label{e5.10}
  \|u-u(x_0)\|_{L^{\infty }(\bar{\Omega }\cap B_r(x_0))}\leq C\left(\|u\|_{L^{\infty }(\Omega)}+\|f\|_{L^n(\Omega)}+\|g\|_{L^{\infty }(\Gamma )}\right)r^{\alpha_0},
\end{equation}
where $0<\alpha_0 <1$ depends only on $n$, $\lambda$, $\Lambda$ and $\delta_0$; $C$ depends also on $\|\gamma \|_{L^{\infty} (\Gamma) }$ and $\check{C}^{-1}$ depends also on the $C^{1}$ modulus of $\Gamma$ at $x_0$.
\end{theorem}
\proof Rewrite the equations as
  \begin{equation*}
    \left\{
    \begin{aligned}
      &u\in S(\lambda,\Lambda,f)~~~~~ \mbox{in}~~ \Omega ; \\
      &\beta  \cdot Du=g-\gamma u~~\mbox{on}~~\Gamma.
    \end{aligned}
    \right.
  \end{equation*}
Then from \Cref{th5.1}, $u$ is $C^{\alpha_0}$ at $x_0$ and\cref{e5.10} holds. \qed~\\

Combining \Cref{th_ca.1} with the interior H\"{o}lder estimate, the boundary local H\"{o}lder estimate (\Cref{th5.2}) follows easily (see \cite[Proposition 4.10]{C-C}).

\section{ Uniqueness and existence of viscosity solutions}\label{S3}
In this section, we derive the Jensen's type uniqueness of viscosity solutions which will be also used to prove the $C^{1,\alpha}$ regularity in a spherical cap in next section. That is, the subsolution minus the supersolution is also a ``subsolution''.  This leads to the uniqueness of viscosity solutions by combining with the A-B-P maximum principle. In addition, we prove an existence result which will be also used to construct auxiliary functions in later sections. The results of this section have been motivated by \cite{C-I-L} and \cite{Is1}.

We start with the following notations (see \cite{C-I-L} or \cite{Is1}).  For $u$ defined on $\bar{\Omega}$ and $x_0\in\bar{\Omega}$, let
\begin{equation*}
\begin{aligned}
   J^{2,+}u(x_0):=&\bigg\{(p,A)\in R^n\times S^n\big|u(x+h)\leq u(x)+p\cdot h+\frac{1}{2}h^{T}Ah+o(|h|^2) \\
  &\mbox{for}~~x+h\in\bar{\Omega}~~\mbox{as}~~h\rightarrow 0\bigg\}
\end{aligned}
\end{equation*}
and
\begin{equation*}
\begin{aligned}
   J^{2,-}u(x_0):=&\bigg\{(p,A)\in R^n\times S^n\big|u(x+h)\geq u(x)+p\cdot h+\frac{1}{2}h^{T}Ah+o(|h|^2) \\
  &\mbox{for}~~x+h\in\bar{\Omega}~~\mbox{as}~~h\rightarrow 0\bigg\}.
\end{aligned}
\end{equation*}
We also define the following:
\begin{equation*}
\begin{aligned}
    \bar{ J}^{2,+}u(x_0):=&\Big\{(p,A)\big|\mbox{ there exist a sequence } (x_m,p_m,A_m) \mbox{ such that }\\
   & (p_m,A_m)\in J^{2,+}u(x_m), x_m\rightarrow x_0, p_m\rightarrow p \mbox{ and } A_m\rightarrow A\Big\}.
\end{aligned}
\end{equation*}
$\bar{ J}^{2,-}u(x_0)$ is defined similarly.

Upon above notations, we have the following observation:
\begin{proposition}\label{l2.4}
Suppose that $\Gamma\in C^2$. Then $u$ is a viscosity subsolution of\cref{e1.1} if and only if
\begin{equation}\label{e3.12}
F(A)\geq f(x_0) ~~~~\forall x_0\in \Omega ,\forall (p,A)\in \bar{J}^{2,+}u(x_0)
\end{equation}
and
\begin{equation}\label{e3.13}
  \left\{
  \begin{aligned}
    &F(A)\geq f(x_0) \\
    &\mbox{or}~~~~~~~~~~~~~~~~~~~~~~~~~~~~~~~ ~~~~~&&,~~\forall x_0\in \Gamma, \forall (p,A)\in \bar{J}^{2,+}u(x_0).\\
  &\beta (x_0)\cdot p+\gamma (x_0) u(x_0)\geq g(x_0)&&
  \end{aligned}
  \right.
\end{equation}
\end{proposition}
\proof  The ``only if'' part is trivial and we only prove the ``if'' part. It is obvious that $F(D^2u (x_0))\geq f(x_0)$ in the viscosity sense if $x_0\in  \Omega$. Hence, we only need to consider the case $x_0\in \Gamma$, let $P$ be a paraboloid touching $u$ by above at $x_0$. We need to prove
\begin{equation*}\label{e3.11}
  \beta(x_0)\cdot DP(x_0)+\gamma (x_0) P(x_0)\geq g(x_0).
\end{equation*}
Suppose not. Without loss of generality, we assume that $x_0=0$ and $\vec{n}(0)=e_n$. Let $\varphi $ denote the representation function of $\Gamma$ near $0$ with $\varphi (0)=0$ and $D\varphi (0)=0$. Then by an argument similar to the proof of \Cref{le2.3}, for any $N>0$, there exist $\eta,\tau '>0$ such that
\begin{equation*}
P (x)\leq P (x',\varphi(x'))+(A-\eta/2)(x_n-\varphi(x'))-N(x_n-\varphi(x'))^2 ~~\forall x\in \bar{\Omega}\cap B_{\tau'}.
\end{equation*}
where $A=\left(g(0)-\gamma(0)P(0)-\beta'(0)\cdot D_{\_}P(0)\right)/\beta _{n}(0)$. Here, we choose $\tau '$ small enough such that $\Omega \cap B_{\tau'}\subset \subset \Omega \cup \Gamma$ and $\beta_n\geq \delta_0/2$ on $\Gamma\cap B_{\tau'}$. For $ \varepsilon > 0$, let $\psi_{\varepsilon } (x)=P (x',\varphi(x'))+(A-\eta/2)(x_n-\varphi(x'))-N(x_n-\varphi(x'))^2+\varepsilon |x|^2$. Then
\begin{equation*}
     \left (D\psi _{\varepsilon}(0), D^2\psi _{\varepsilon}(0)\right)\in J^{2,+}u(0).
\end{equation*}
Hence, we have
\begin{equation*}
  F(D^2\psi _{\varepsilon}(0))\geq f(0)
\end{equation*}
or
\begin{equation*}
  \beta (0)\cdot D\psi _{\varepsilon}(0)+\gamma (0) \psi _{\varepsilon}(0)\geq g(0).
\end{equation*}
Similar to the argument in the proof of \Cref{le2.3}, by choosing $N$ large enough and $\varepsilon$ small enough, we obtain a contradiction. \qed~\\

\begin{remark}\label{re3.1}
In some context, that\cref{e3.12} and\cref{e3.13} hold is the definition of viscosity subsolution (cf. \cite[Definition 7.4]{C-I-L}). This proposition states that these two definitions are equivalent if $\Gamma\in C^2$.
\end{remark}

We introduce the following two results (see \cite[Theorem 4.1]{Is1} and \cite[Lemma 3.1]{Is1}).
\begin{proposition}\label{ishii}
Let $\Omega$ be a bounded domain, $\partial\Omega\supset \Gamma\in C^2$ be relatively open and $\beta\in C^2(\bar{\Gamma})$ with $\beta\cdot \vec{n}\geq \delta_0$ on $\Gamma$. Given $x_0\in\Gamma$, there are positive constants $r_0$, $C$ and a family of $\{w_{\varepsilon }\}_{\varepsilon > 0}$ of $C^{1,1}$ functions on $B_{r_0}(x_0)\times B_{r_0}(x_0)$ such that for $\varepsilon > 0$ and $x,y\in B_{r_0}(x_0)$,
\begin{equation}\label{ish}
\begin{aligned}
 w_{\varepsilon }(x,x)\leq &\varepsilon,  \\
 w_{\varepsilon }(x,y)\geq & \frac{|x-y|^2}{8\varepsilon },\\
 \beta (x) \cdot D_{x}w_{\varepsilon }(x,y)\leq& C\left (\frac{|x-y|^2}{\varepsilon }+\varepsilon \right)~~\mbox{if}~~x\in \Gamma~~\mbox{and}~~y\in \bar{\Omega},\\
 \beta (y) \cdot D_{y}w_{\varepsilon }(x,y)\leq & C\left (\frac{|x-y|^2}{\varepsilon }+\varepsilon \right)~~\mbox{if}~~y\in \Gamma~~\mbox{and}~~x\in \bar{\Omega},\\
\end{aligned}
\end{equation}
and
\begin{equation}\label{e3.4}
  \left(Dw_{\varepsilon }(x,y),\frac{C}{\varepsilon }
 \left(
 \begin{array}{cc}
 I & -I \\
 -I & I \\
 \end{array}
 \right)
 +C\left(\frac{|x-y|^2}{\varepsilon }+\varepsilon\right)
  \left(
 \begin{array}{cc}
 I & 0 \\
 0 & I \\
\end{array}
\right)
 \right)\in J^{2,+}w_{\varepsilon }(x,y).
\end{equation}
\end{proposition}

\begin{proposition}\label{ishii2}
Let $u\in USC(\bar{\Omega })$, $v\in LSC(\bar{\Omega })$ and $w(x,y)=u(x)-v(y)$ for $x,y\in\bar{\Omega }$. Assume that
\begin{equation*}
\left((p,q),\tau
 \left(
 \begin{array}{cc}
 I & -I \\
 -I & I \\
 \end{array}
 \right)
 +\nu
  \left(
 \begin{array}{cc}
 I & 0 \\
 0 & I \\
\end{array}
\right)
 \right)\in J^{2,+}w(\hat{x},\hat{y})
\end{equation*}
for some $p,q\in R^n $, $\hat{x},\hat{y}\in \bar{\Omega }$, $\tau > 0$ and $\nu >0$. Then there are $X,Y\in S^{n}$ for which
\begin{equation*}
  -3\tau \left(
           \begin{array}{cc}
             I & 0 \\
             0 & I \\
           \end{array}
         \right)
\leq
\left(
  \begin{array}{cc}
    X-\nu I & 0 \\
    0 & Y-\nu I \\
  \end{array}
\right)
\leq 3\tau
\left(
  \begin{array}{cc}
    I & -I \\
    -I & I \\
  \end{array}
\right)
\end{equation*}
and
\begin{equation*}
  (p,X)\in \bar{J}^{2,+}u(\hat{x})~~\mbox{and}~~(-q,-Y)\in \bar{J}^{2,-}v(\hat{y}).
\end{equation*}
\end{proposition}

Combining above two results and the doubling variable arguments (see \cite[Section 3]{C-I-L}), we prove a Jensen's type uniqueness result.
\begin{theorem}\label{th3.11}
Suppose that $\Gamma\in C^2$ and $\beta\in C^2(\bar{\Gamma})$. Let $u$ and $v$ satisfy
\begin{equation*}
\left\{\begin{aligned}
 &F(D^2u)\geq f_1 &&\mbox{in}~~\Omega; \\
 &\beta\cdot Du+\gamma u\geq g_1 &&\mbox{on}~~\Gamma
\end{aligned}
\right.
\end{equation*}
and
\begin{equation*}
\left\{\begin{aligned}
 &F(D^2v)\leq f_2 &&\mbox{in}~~\Omega; \\
 &\beta\cdot Dv+\gamma v\leq g_2 &&\mbox{on}~~\Gamma.
\end{aligned}
\right.
\end{equation*}

Then
    \begin{equation*}
    \left\{
    \begin{aligned}
     &u-v\in \underline{S } (\lambda/n,\Lambda,f_1-f_2 ) &&\mbox{in}~~\Omega; \\
     &\beta\cdot D(u-v)+\gamma (u-v)\geq g_1-g_2 &&\mbox{on}~~\Gamma.
    \end{aligned}
    \right.
  \end{equation*}
\end{theorem}
\proof
Let $P$ be a paraboloid touching $u-v$ by above at $x_0\in  \Omega$. We need to prove
\begin{equation}\label{e3.5}
M^+(D^2P, \lambda /n,\Lambda)\geq f_1(x_0)-f_2(x_0).
\end{equation}

Let $P_{\varepsilon _{0}}=P+\varepsilon _{0}|x-x_0|^2$, $w(x,y)=u(x)-v(y)-P_{\varepsilon _{0}}(y)$, $\varphi_{\alpha } (x,y)=\alpha |x-y|^2/2$ ($\alpha>0$) and $(x_{\alpha },y_{\alpha })$ be a maximum point of $w-\varphi_{\alpha } $. Then
\begin{equation*}
\left(\Big(\alpha (x_{\alpha }- y_{\alpha }),\alpha (y_{\alpha }- x_{\alpha })\Big),\alpha
 \left(
 \begin{array}{cc}
 I & -I \\
 -I & I \\
 \end{array}
 \right)
 \right)\in J^{2,+}w(x_{\alpha },y_{\alpha }).
\end{equation*}
Applying \Cref{ishii2} with $\tau=\alpha$ and $\nu=0$, we conclude that there exist $X_{\alpha}$, $Y_{\alpha}\in S^n $ such that
\begin{equation*}
  \left(\alpha (x_{\alpha }- y_{\alpha} ),X_{\alpha }\right)\in \bar{J}^{2,+}u(x_{\alpha }),~~\left(\alpha (x_{\alpha }- y_{\alpha} ),Y_{\alpha }\right)\in \bar{J}^{2,-}(v+P_{\varepsilon _{0}})(y_{\alpha })
\end{equation*}
and
\begin{equation}\label{e3.2}
\left(
  \begin{array}{cc}
    X_{\alpha } & 0 \\
    0 & -Y_{\alpha } \\
  \end{array}
\right)
\leq 3\alpha
\left(
  \begin{array}{cc}
    I & -I \\
    -I & I \\
  \end{array}
\right).
\end{equation}

Recall that $x_0$ is the unique maximum point of $u-v-P_{\varepsilon _0}$. Then it is easy to verify that $x_{\alpha }\rightarrow x_0$ and $y_{\alpha }\rightarrow x_0$ as $\alpha \rightarrow + \infty $ (see \cite[Lemma 3.1]{C-I-L}). This implies that $x_{\alpha },y_{\alpha }\in  \Omega$ for $\alpha$ large enough.  Since $F(D^2u)\geq f_1$ and $F(D^2v)\leq f_2$, we have that $F(X_{\alpha })\geq f_1(x_{\alpha })$ and $F(Y_{\alpha }-D^2P-2\varepsilon _0 I)\leq f_2(y_{\alpha })$. Note that\cref{e3.2} implies $X_{\alpha }\leq Y_{\alpha }$. Hence,
\begin{equation*}
  \begin{aligned}
    f_2(y_{\alpha }) &\geq  F(Y_{\alpha }-D^2P-2\varepsilon _0 I) \\
     & \geq F(Y_{\alpha })-M^+(D^2P, \lambda /n,\Lambda)-2n\Lambda\varepsilon _{0}\\
     & \geq F(X_{\alpha })-M^+(D^2P, \lambda /n,\Lambda)-2n\Lambda\varepsilon _{0}\\
     & \geq f_1(x_{\alpha})-M^+(D^2P, \lambda /n,\Lambda)-2n\Lambda\varepsilon _{0}
  \end{aligned}
\end{equation*}
Let $\alpha \rightarrow \infty $ and it follows that
\begin{equation*}
  M^+(D^2P, \lambda /n,\Lambda)\geq f_1(x_0)-f_2(x_0)-2n\Lambda\varepsilon _{0}.
\end{equation*}
Next, let $\varepsilon _{0}\rightarrow 0$ and\cref{e3.5} follows.

In the following, we consider the boundary case. Let $P$ be a paraboloid touching $u-v$ by above at $x_0\in \Gamma$. Without loss of generality, we assume that $x_0=0$ and $\vec{n}(0)=e_n$. Let $\varphi $ denote the representation function of $\Gamma$ near $0$ with $\varphi (0)=0$ and $D\varphi (0)=0$. We need to prove
\begin{equation}\label{ze1.1}
  \beta (0) \cdot DP(0)+\gamma(0)P(0)\geq g_1(0)-g_2(0).
\end{equation}
Suppose not. Then
\begin{equation*}
  P_n (0)<\left(g_1(0)-g_2(0)-\gamma(0) P (0)-\beta'(0)\cdot D_{\_}P (0)\right)/\beta_n(0):=A_1-A_2,
\end{equation*}
where $A_1=(g_1(0)-\gamma(0)u(0))/\beta_n(0)$ and $A_2=(g_2(0)-\gamma(0)v(0)+\beta'(0)\cdot D_{\_}P(0))/\beta_n(0)$.
By the continuity of $P_n $, there exist $ \tau > 0$ and $\eta> 0$ such that
\begin{equation*}
 P_n (x)<A_1-A_2-3\eta~~\forall x\in \Gamma\cap B_\tau.
\end{equation*}
By Taylor's formula, for any $x\in\bar{\Omega}\cap B_\tau$,
\begin{equation*}
\begin{aligned}
    P (x)=&P (x',\varphi(x'))+P _{n}(x',\varphi(x'))(x_n-\varphi(x'))+\frac{1}{2}P _{nn}(x_n-\varphi(x'))^2\\
\leq &P (x',\varphi(x'))+(A_1-A_2-3\eta)(x_n-\varphi(x'))+\frac{1}{2}P_{nn}(x_n-\varphi(x'))^2.
\end{aligned}
\end{equation*}
By the boundedness of $P_{nn} $, for any $N>0$ there exists $\tau '> 0$ such that
\begin{equation*}
P (x)\leq P (x',\varphi(x'))+(A_1-A_2-2\eta)(x_n-\varphi(x'))-2N(x_n-\varphi(x'))^2~~\forall x\in \bar{\Omega}\cap B_{\tau'}.
\end{equation*}
The constant $N$ is large enough to be chosen later. We choose $\tau '<r_0$ small enough such that $\bar{\Omega}\cap B_{\tau'}\subset \subset \Omega\cup \Gamma$ and $\beta_n\geq \delta_0/2$ on $\Gamma\cap B_{\tau'}$ where $r_0$ is as in \Cref{ishii}.

Let $\psi(x)=(A_1-\eta)(x_n-\varphi(x'))-N(x_n-\varphi(x'))^2$. For $\varepsilon_1 > 0$, let $\psi_{\varepsilon_1 } (x)=P (x',\varphi(x'))-(A_2+\eta)(x_n-\varphi(x'))-N(x_n-\varphi(x'))^2+\varepsilon_1 |x|^2$. Then
\begin{equation*}
  \psi+\psi_{\varepsilon_1 } > P \geq u-v~~\mbox{in}~~\Omega\cap B_{\tau'}\backslash \{0\}~~\mbox{and}~~\psi+\psi_{\varepsilon_1 } =P =u-v \mbox{ at } 0.
\end{equation*}
Next, for $\theta >0$, let
\begin{equation*}
  \tilde{u}(x)=u(x)+\theta ^2x_n-\theta |x|^2-\psi(x)
\end{equation*}
and
\begin{equation*}
  \tilde{v}(x)=v(x)-\theta ^2x_n+\psi _{\varepsilon _1}(x).
\end{equation*}
Finally, for $\varepsilon >0$, let
\begin{equation*}
  \Phi_{\varepsilon } (x,y)=\tilde{u}(x)-\tilde{v}(y)-w_{\varepsilon }(x,y),
\end{equation*}
where $w_{\varepsilon }(x,y)$ is as in \Cref{ishii}. Then by\cref{ish},
\begin{equation*}
 \begin{aligned}
 -\varepsilon \leq& \Phi_{\varepsilon } (0,0)\leq \Phi_{\varepsilon } (x_\varepsilon ,y_\varepsilon )\\
 =&u(x_{\varepsilon })+\theta^2x_{\varepsilon ,n}-\theta|x_{\varepsilon }|^2-\psi(x_{\varepsilon })-\Big(v(y_{\varepsilon })-\theta^2y_{\varepsilon ,n}+\psi_{\varepsilon _1}(y_{\varepsilon })\Big)-w_{\varepsilon }(x_{\varepsilon },y_{\varepsilon })\\
 \leq & u(x_{\varepsilon })-v(y_{\varepsilon })-\theta |x_{\varepsilon }-\theta e_n|^2+\theta ^{3}-\theta^2  (x_{\varepsilon }-y_{\varepsilon })\cdot e_n\\
 &- \frac{|x_{\varepsilon }-y_{\varepsilon }|^2}{8 \varepsilon }-\psi(x_{\varepsilon })-\psi _{\varepsilon _1}(y_{\varepsilon}),
 \end{aligned}
\end{equation*}
 where $(x_{\varepsilon },y_{\varepsilon })$ is a maximum point of $\Phi_{\varepsilon }$. It follows that for $\varepsilon _1$ and $\theta $ fixed,
 \begin{equation*}
  |x_{\varepsilon }-y_{\varepsilon }|\longrightarrow 0~~\mbox{as}~~\varepsilon \rightarrow 0
 \end{equation*}
and
\begin{equation}\label{e3.8}
 \underset{\varepsilon \rightarrow 0} {\overline{ \lim}}~\frac{|x_{\varepsilon }-y_{\varepsilon }|^2}{8\varepsilon }\leq \theta ^{3},~~x_{\varepsilon },y_{\varepsilon }\rightarrow B_{\theta}(\theta e_n)~~\mbox{as}~~\varepsilon \rightarrow 0.
\end{equation}
We choose $\varepsilon $ and $\theta $ small enough such that $x_{\varepsilon },y_{\varepsilon }\in \bar{\Omega}\cap B_{\tau'}$.

Since $(x_{\varepsilon },y_{\varepsilon })$ is a maximum point of $\Phi_{\varepsilon }$, we have
\begin{equation*}
  \tilde{u}(x)-\tilde{v}(y)\leq w_{\varepsilon }(x,y)+\tilde{u}(x_{\varepsilon })-\tilde{v}(y_{\varepsilon })-w_{\varepsilon }(x_{\varepsilon },y_{\varepsilon })~~\forall x,y\in  \bar{\Omega}\cap B_{\tau'}.
\end{equation*}
Using\cref{e3.4} in \Cref{ishii}, we have
\begin{equation*}
   \left((p,q),\frac{C}{\varepsilon }
 \left(
 \begin{array}{cc}
 I & -I \\
 -I & I \\
 \end{array}
 \right)
 +Cs
  \left(
 \begin{array}{cc}
 I & 0 \\
 0 & I \\
\end{array}
\right)
 \right)\in J^{2,+}\tilde{w}(x_{\varepsilon },y_{\varepsilon }) ,
\end{equation*}
where $\tilde{w}(x,y)=\tilde{u}(x)-\tilde{v}(y)$, $p=D_{x}w_{\varepsilon }(x_{\varepsilon },y_{\varepsilon })$, $q=D_{y}w_{\varepsilon }(x_{\varepsilon },y_{\varepsilon })$  and $s=\left (|x_{\varepsilon }-y_{\varepsilon }|^2/\varepsilon +\varepsilon \right)$. By \Cref{ishii2}, there are $X,Y\in S^{n}$ such that
\begin{equation}\label{e1.3}
  \left(
  \begin{array}{cc}
    X-Cs I & 0 \\
    0 & Y-Cs I \\
  \end{array}
\right)
\leq \frac{3C}{\varepsilon }
\left(
  \begin{array}{cc}
    I & -I \\
    -I & I \\
  \end{array}
\right)
\end{equation}
and
\begin{equation*}
  (p,X)\in \bar{J}^{2,+}\tilde{u}(x_{\varepsilon })~~\mbox{and}~~(-q,-Y)\in \bar{J}^{2,-}\tilde{v}(y_{\varepsilon }).
\end{equation*}
By the definition of $\tilde{u}$, $\tilde{v}$, $\psi$ and $\psi _{\varepsilon _1}$, we have
\begin{equation}\label{e1.4}
\left(\tilde {p},X+2\theta I+\bar {M}\right)\in \bar{J}^{2,+}u(x_{\varepsilon })
\end{equation}
and
\begin{equation}\label{e1.2}
\begin{aligned}
 &\left(\tilde {q}, -Y-\tilde {M}-2\varepsilon _1 I\right)\in \bar{J}^{2,-}v(y_{\varepsilon })
\end{aligned}
\end{equation}
where
\begin{equation*}
  \tilde {p}=p-\theta ^{2} e_n+ 2\theta x_{\varepsilon }+\Big((A_1-\eta)-2N\big(x_{\varepsilon, n}-\varphi (x_{\varepsilon }')\big)\Big)(-D\varphi (x_{\varepsilon }'),1),
\end{equation*}
\begin{equation*}
\begin{aligned}
  \tilde {q}=&-q+\theta ^{2} e_n+\Big((A_2+\eta)+2N\big(y_{\varepsilon,n}-\varphi (y_{\varepsilon}')\big)\Big)(-D\varphi (y_{\varepsilon}'),1)- 2\varepsilon_1 y_{\varepsilon }\\
&-\Big(\big(D_{\_}P(y_{\varepsilon}' ,\varphi (y_{\varepsilon }'))\big)+P_n(y_{\varepsilon}' ,\varphi (y_{\varepsilon }'))D\varphi (y_{\varepsilon }'),0\Big),
\end{aligned}
\end{equation*}
\begin{equation*}
\bar {M}_{ij}=\left\{
  \begin{aligned}
  &(-A_1+\eta-2N\varphi +2Nx_{\varepsilon ,n})\varphi_{ij}-2N\varphi _{i}\varphi _{j},  && i,j<n,\\
  &2N\varphi _{i}, && i<n,j=n,\\
  &-2N, && i=j=n,
  \end{aligned}\right.
\end{equation*}
and
\begin{equation*}
\tilde {M}_{ij}=\left\{
  \begin{aligned}
  &(A_2+\eta-2N\varphi+2Ny_{\varepsilon ,n})\varphi_{ij}-2N\varphi _{i}\varphi _{j} \\
  &+P_{ij}+P_{in}\varphi _{j}+P_{jn}\varphi _{i}+P_n\varphi_{ij}, && i,j<n,\\
  &2N\varphi _{i}, && i<n,j=n,\\
  &-2N, && i=j=n.
  \end{aligned}\right.
\end{equation*}

If $x_{\varepsilon }\in \Gamma$, by\cref{e1.4} and the definition of viscosity solutions,
\begin{equation}\label{e3.9}
  \beta (x_{\varepsilon})\cdot\tilde{p}+\gamma(x_{\varepsilon })u(x_{\varepsilon })\geq g_1(x_{\varepsilon }),
\end{equation}
i.e.,
\begin{equation*}
\begin{aligned}
  &\beta (x_{\varepsilon})\cdot p-\theta^2\beta_n(x_{\varepsilon })+2\theta \beta(x_{\varepsilon })\cdot x_{\varepsilon }\\
+&\Big((A_1-\eta)-2N\big(x_{\varepsilon, n}-\varphi (x_{\varepsilon }')\big)\Big)(-\beta'(x_{\varepsilon})\cdot D\varphi (x_{\varepsilon }')+\beta_n(x_{\varepsilon})) \\
+&\gamma(x_{\varepsilon })u(x_{\varepsilon })\geq g_1(x_{\varepsilon }).
\end{aligned}
\end{equation*}
By\cref{ish} and\cref{e3.8},
\begin{equation*}
\beta (x_{\varepsilon }) \cdot p \leq C\left (\frac{|x_{\varepsilon }-y_{\varepsilon }|^2}{\varepsilon }+\varepsilon \right)\leq C(\theta^3+\varepsilon)
\end{equation*}
if $\varepsilon$ and $\theta$ are small. On the other hand, since $\varphi (0)=0$, $D\varphi (0)=0$, $\beta_n(0)\geq \delta_0$ and $x_{\varepsilon }\rightarrow 0$ as $\varepsilon ,\theta \rightarrow 0$ (recall\cref{e3.8}), by choosing $\varepsilon $ and $\theta $ small enough, we have (for $N$ fixed)
\begin{equation*}
(A_1-\eta/2)\beta_n(x_{\varepsilon}) +\gamma(x_{\varepsilon })u(x_{\varepsilon })\geq g_1(x_{\varepsilon }).
\end{equation*}
Since $A_1=(g_1(0)-\gamma(0)u(0))/\beta_n(0)$, by choosing $\varepsilon $ and $\theta $ small enough and the continuity of $u$, $\beta$, $\gamma$ and $g_1$, we have
\begin{equation*}
-\eta\delta_0/4\geq 0,
\end{equation*}
which is a contradiction. Thus, $x_{\varepsilon }\in \Omega$ and
\begin{equation}\label{e3.6}
    F(X+2\theta I+\bar {M})\geq f_1(x_{\varepsilon }).
\end{equation}

By a similar argument for $y_{\varepsilon }$, we have that $y_{\varepsilon }\in \Omega$ and
\begin{equation}\label{e3.7}
  F(-Y-\tilde{M}- 2\varepsilon _1I)\leq f_2(y_{\varepsilon }).
\end{equation}
Note that\cref{e1.3} implies,
\begin{equation*}
  X-CsI\leq -Y+CsI,~~\mbox{i.e.,}~~X -2CsI\leq -Y.
\end{equation*}
Then from\cref{e3.6} and\cref{e3.7},
\begin{equation*}
\begin{aligned}
  f_2(y_{\varepsilon }) & \geq F(-Y-\tilde{M}- 2\varepsilon _1I) \\
   & \geq F(-Y)-M^+(\tilde{M},\lambda/n,\Lambda)-2n\Lambda \varepsilon _1\\
   & \geq F(X-2CsI)-M^+(\tilde{M},\lambda/n,\Lambda)-2n\Lambda \varepsilon _1\\
   & \geq F(X)-M^+(\tilde{M},\lambda/n,\Lambda)-2n\Lambda \varepsilon _1-2n\Lambda Cs\\
   & \geq f_1(x_{\varepsilon })-M^+(\bar{M},\lambda/n,\Lambda)-M^+(\tilde{M},\lambda/n,\Lambda)-2n\Lambda(\varepsilon _1+Cs+\theta).
\end{aligned}
\end{equation*}
Recall the definitions of $\bar{M}$ and $\tilde{M}$. By choosing $N$ large enough (independent of $\theta$ and $\varepsilon $), and $\theta$ and $\varepsilon $ small enough, we obtain a contradiction (note that $\varphi (0)=0$, $D\varphi (0)=0$ and $x_{\varepsilon },y_{\varepsilon }\rightarrow 0$ as $\varepsilon, \theta \rightarrow 0 $). Therefore,\cref{ze1.1} holds.\qed~\\

\begin{remark}\label{re.remove}
If $u$ or $v$ belongs to $C^1$, the conditions $\Gamma\in C^2$ and $\beta\in C^2(\Gamma)$ can be removed; and $\beta\cdot D(u-v)+\gamma (u-v)\geq g_1-g_2$ on $\Gamma$ can be verified directly by the definition of viscosity solutions.
\end{remark}

Combining above theorem with the A-B-P maximum principle, we derive a uniqueness result:
\begin{theorem}\label{th.uni}
Let $\Gamma\in C^2$, $\beta\in C^2(\bar{\Gamma})$, $\gamma\leq 0$ and $\varphi \in C(\partial \Omega \backslash \Gamma) $. Suppose that there exists $\xi\in\partial B_1$ such that
\begin{equation*}\label{e.exis.2}
  \beta\cdot \xi \geq \delta_1 ~\mbox{on}~\Gamma.
\end{equation*}

Then there exists at most one viscosity solution of
  \begin{equation*}
    \left\{
    \begin{aligned}
      &F(D^2u)=f &&\mbox{in}~~ \Omega; \\
      &\beta  \cdot Du +\gamma u = g &&\mbox{on}~~\Gamma;\\
      &u=\varphi  &&\mbox{on}~~\partial\Omega\backslash\Gamma.
    \end{aligned}
    \right.
  \end{equation*}
\end{theorem}

Next, we use Perron's method to prove an existence result for fully nonlinear elliptic equations with a ``small'' portion of oblique boundary.
\begin{theorem}\label{th.existence}
Let $\Gamma\in C^2$, $\gamma\leq 0$ and $\varphi \in C(\partial \Omega \backslash \Gamma)$. Suppose that there exists $\xi\in\partial B_1$ such that
\begin{equation}\label{existence_4}
  \beta\cdot \xi \geq \delta_1 ~\mbox{on}~\Gamma.
\end{equation}
In addition, suppose that $\Omega$ satisfies an exterior cone condition at any $x\in \partial\Omega\backslash \bar{\Gamma}$ and satisfies an exterior sphere condition at any $x\in \bar{\Gamma}\cap (\partial\Omega\backslash \Gamma)$.

Then there exists a unique viscosity solution $u\in C(\bar{\Omega })$ of
  \begin{equation}\label{existence}
    \left\{
    \begin{aligned}
      &F(D^2u)=f &&\mbox{in}~~ \Omega; \\
      &\beta  \cdot Du+\gamma u = g &&\mbox{on}~~\Gamma;\\
      &u=\varphi  &&\mbox{on}~~\partial\Omega\backslash\Gamma.
    \end{aligned}
    \right.
  \end{equation}
\end{theorem}

\proof We assume that $F(0)=0$. Otherwise, by the uniform ellipticity, there exists $t\in R$ such that $F(t\delta_{nn})=0$ and $|t|\leq  |F(0)|/\lambda$. Let $G(M):=F(M+t\delta_{nn})$ and then $G(0)=0$. Hence, $\tilde{u}+tx_n^2/2$ is the unique solution of\cref{existence} where $\tilde{u}$ is the unique solution of
\begin{equation*}
    \left\{
    \begin{aligned}
      &G(D^2\tilde{u})=f &&\mbox{in}~~\Omega ; \\
      &\beta \cdot D\tilde{u}+\gamma \tilde{u} =g-tx_n\beta_n-tx_n^2\gamma &&\mbox{on}~~\Gamma;\\
      &\tilde{u}=\varphi-tx_n^2/2 &&\mbox{on}~~\partial \Omega \backslash \Gamma.
    \end{aligned}
    \right.
\end{equation*}

From now on, we assume that $F(0)=0$. Let
\begin{equation*}
  \textbf{V}:=\left\{v\in USC(\bar{\Omega })\big| v \mbox{ is a subsolution of\cref{existence} with } v\leq \varphi \mbox{ on }\partial \Omega \backslash \Gamma \right\}.
\end{equation*}
By choosing proper positive constants $K_1$, $K_2$ and $K_3$, $K_1x_n^2+K_2x_n-K_3\in \textbf{V}$ and hence $\textbf{V}$ is nonempty. Set
\begin{equation*}
  w(x)=\sup_{v\in\textbf{V}}v(x),w^{\ast}(x)=\underset{y\rightarrow x}{\overline{\lim}}w(y) \mbox{ and } w_{\ast}(x)=\underset{y\rightarrow x}{\underline{\lim}}w(y).
\end{equation*}
Then, $w^{\ast}\in USC(\bar{\Omega })$, $w_{\ast}\in LSC(\bar{\Omega })$ and $w_{\ast}\leq w\leq w^{\ast}$. Next, we divide the proof into three steps.~\\

\emph{Step 1}. We prove that $w^{\ast}$ is a subsolution.

For any $ \tilde{x}\in \Omega $ and paraboloid $P$ touching $w^{\ast}$ by above at $\tilde{x}$, we need to prove that $F(D^2P)\geq f(\tilde{x})$. Suppose not. Then there exists $\varepsilon > 0$ such that
\begin{equation}\label{e2.9}
  F(D^2P+2\varepsilon I)< f(\tilde{x})-2\varepsilon .
\end{equation}
By the definition of $w^{\ast}$, there exist $\{x_k\}\subset \Omega $ and $\{v_k\}\subset \textbf{V}$ such that
\begin{equation*}
 x_k\rightarrow \tilde{x} \mbox{ and } w^{\ast}(\tilde{x})=\lim_{k\rightarrow \infty} v_k(x_k).
\end{equation*}
Then there exists $ r>0$, for $k$ large enough, such that $x_k\in B_{r/2}(\tilde{x})\subset B_{r}(\tilde{x})\subset \subset \Omega$, $|f(x)-f(\tilde{x})|\leq \varepsilon $ for any $x\in B_{r}(\tilde{x})$ and
\begin{equation*}
|P(x_k)+\varepsilon |x_k-\tilde{x}|^2-v_k(x_k)|\leq|P(x_k)-P(\tilde{x})|+\varepsilon |x_k-\tilde{x}|^2+|P(\tilde{x})-v_k(x_k)|< \frac{1}{2}\varepsilon r^2.
\end{equation*}
On the other hand, $P+\varepsilon |x-\tilde{x}|^2 -v_k\geq \varepsilon r^2~~\mbox{on}~~\partial B_r(\tilde{x})$ (note that $v_k\leq w\leq w^{\ast}$). Hence, $P+\varepsilon |x-\tilde{x}|^2+c_0$ touches $v_k$ by above at some $x^{\ast}\in B_r(\tilde{x})$ for a proper constant $c_0$. Since $v_k$ is a subsolution,
\begin{equation*}
  F(D^2P+2\varepsilon I)\geq f(x^{\ast })\geq f(\tilde{x})-\varepsilon ,
\end{equation*}
which contradicts with\cref{e2.9}.

Next, we consider the case $\tilde{x}\in \Gamma$. Without loss of generality, we assume that $\tilde{x}=0$ and $\vec{n}(0)=e_n$. Let $\varphi $ denote the representation function of $\Gamma$ near $0$ with $\varphi (0)=0$ and $D\varphi (0)=0$. We need to prove that $\beta (0)\cdot DP(0)+\gamma(0)P(0)\geq g(0)$. Suppose not. Then
\begin{equation}\label{e1.36}
  P_n(0)<\frac{1}{\beta _{n}(0)}\left(g(0)-\gamma(0)P(0)-\beta'(0)\cdot D_{\_}P(0)\right):=A.
\end{equation}
By the continuity of $P_n$, there exist $\tau,\eta>0$ such that
\begin{equation*}
  P_n(x)<A-\eta~~\forall x\in B_{\tau}.
\end{equation*}
By Taylor's formula, for any $x\in\bar{\Omega}\cap B_\tau$,
\begin{equation*}
\begin{aligned}
    P (x)=&P (x',\varphi(x'))+P _{n}(x',\varphi(x'))(x_n-\varphi(x'))+\frac{1}{2}P _{nn}(x_n-\varphi(x'))^2\\
\leq &P (x',\varphi(x'))+(A-\eta)(x_n-\varphi(x'))+\frac{1}{2}P_{nn}(x_n-\varphi(x'))^2.
\end{aligned}
\end{equation*}
By the boundedness of $P_{nn} $, for any $N>0$ there exists $\tau '> 0$ such that
\begin{equation*}
P (x)\leq P (x',\varphi(x'))+(A-\eta/2)(x_n-\varphi(x'))-N(x_n-\varphi(x'))^2~~\forall x\in \bar{\Omega}\cap B_{\tau'}.
\end{equation*}
The constant $N$ is large enough to be chosen later. Here, we choose $\tau '$ small enough such that $\Omega \cap B_{\tau'}\subset \subset \Omega \cup \Gamma$ and $\beta_n\geq \delta_0/2$ on $\Gamma\cap B_{\tau'}$. For $ \varepsilon > 0$, let $\psi_{\varepsilon } (x)=P (x',\varphi(x'))+(A-\eta/2)(x_n-\varphi(x'))-N(x_n-\varphi(x'))^2+\varepsilon |x|^2$. Then
\begin{equation*}
  \psi_{\varepsilon } > P \geq w^{\ast} ~~\mbox{in}~~\Omega \cap B_{\tau'}\mbox{ and }\psi_{\varepsilon } (0)=P (0)=w^{\ast }(0).
\end{equation*}
Similar to the interior case, there exist $\{x_k\}\subset \Omega \cup \Gamma$ and $\{v_k\}\subset \textbf{V}$ such that
\begin{equation*}
  x_k\rightarrow 0 \mbox{ and }  w^{\ast}(0)=\lim_{k\rightarrow \infty} v_k(x_k).
\end{equation*}
Then for $k$ large enough,
\begin{equation*}
  x_k\in \bar{\Omega }\cap B_{\varepsilon \tau'/2} \mbox{ and } \psi_{\varepsilon } -v_k\geq \varepsilon^3(\tau ')^2~~\mbox{on}~~\partial B_{\varepsilon \tau'/2}\cap \Omega
\end{equation*}
and
\begin{equation*}
  |v_k(x_k)-\psi_{\varepsilon } (x_k)|\leq |v_k(x_k)-\psi_{\varepsilon } (0)|+|\psi_{\varepsilon } (0)-\psi_{\varepsilon } (x_k)|<\frac{1}{2}\varepsilon^3( \tau ')^{2}.
\end{equation*}
Then $\psi_{\varepsilon }+c_0$ touches $v_k$ by above at some $x^{\ast}\in \bar{\Omega }\cap B_{\varepsilon \tau'}$ for a proper $c_0$. If $x^{\ast}\in \Omega$, then
\begin{equation*}
  F(D^2\psi_{\varepsilon } (x^{\ast }))\geq f(x^{\ast })
\end{equation*}
which is impossible by taking $N$ large enough as in \Cref{th3.11}. If $x^{\ast}\in \Gamma$, then
\begin{equation}\label{existence_2}
  \beta (x^{\ast })\cdot D\psi (x^{\ast })+\gamma(x^{\ast})\psi (x^{\ast })\geq g(x^{\ast }).
\end{equation}
As in \Cref{th3.11}, by recalling the definition of $A$ (see\cref{e1.36}) and the continuity of the functions in\cref{existence_2}, we obtain a contradiction.

From above arguments, we conclude that $w^{\ast }$ is a subsolution. Hence, $w^{\ast }\in \textbf{V}$ and it follows that $w^{\ast }\leq w$. Recall that $w^{\ast }\geq w$. Hence, $w=w^{\ast }$ is a subsolution.~\\

\emph{Step 2}. We prove that $w_{\ast }$ is a supersolution.

Suppose not. Then there exist $\tilde{x}\in \Omega \cup \Gamma$ and a paraboloid $P$ touching $w_{\ast }$ by below at $\tilde{x}$ such that
\begin{equation*}
  F(D^2P)>f(\tilde{x})~~\mbox{if}~~\tilde{x}\in \Omega
\end{equation*}
and
\begin{equation*}
  \beta (\tilde{x})\cdot DP(\tilde{x})+\gamma(\tilde{x})P(\tilde{x})>g(\tilde{x})~~\mbox{if}~~\tilde{x}\in \Gamma.
\end{equation*}

If $\tilde{x}\in \Omega $, take $\varepsilon _0>0$ small such that
\begin{equation}\label{e1.37}
  F(D^2P-2\varepsilon _0I)>f(\tilde{x}).
\end{equation}
Take $r>0$ small such that $B_r(\tilde{x})\subset\subset \Omega$. Let
\begin{equation*}
  \psi _{\varepsilon _0}=P-\varepsilon _0|x-\tilde{x}|^2+\frac{1}{2}\varepsilon _0r^2.
\end{equation*}
Then
\begin{equation*}
w\geq w_{\ast }\geq \psi _{\varepsilon _0}~~\mbox{on}~~\partial B_r(\tilde{x}) \mbox{ and }  w_{\ast }(\tilde{x})<\psi _{\varepsilon _0}(\tilde{x}),
\end{equation*}
which implies that there exists $x_1\in B_r(\tilde{x})$ such that
\begin{equation}\label{e1.38}
  w(x_1)<\psi _{\varepsilon _0}(x_1).
\end{equation}
Define
\begin{equation*}
  \bar{w}=\left\{
  \begin{aligned}
    \max & (w,\psi _{\varepsilon _0})  &&\mbox{if}~~x\in \bar{B}_r(\tilde{x});\\
    w &  &&\mbox{if}~~x\notin \bar{B}_r(\tilde{x}).
  \end{aligned}
  \right.
\end{equation*}
It is easy to verify that $\bar{w}$ is a subsolution (recall that $w$ is a subsolution). Hence, $\bar{w}\leq w$. In particular, $\psi _{\varepsilon _0}(x_1)\leq \bar{w}(x_1)\leq w(x_1)$ which contradicts with\cref{e1.38}.

If $\tilde{x}\in \Gamma$, similar to previous arguments, we assume that $\tilde{x}=0$ and $\vec{n}(0)=e_n$. Then, there exist $\eta,\tau '> 0$ such that
\begin{equation*}
P (x)> P (x',\varphi(x'))+(A+\eta/2)(x_n-\varphi(x'))+N(x_n-\varphi(x'))^2~~\forall x\in \bar{\Omega}\cap B_{\tau'}.
\end{equation*}
where $A:=\left(g(0)-\gamma(0)P(0)-\beta' (0)\cdot D_{\_}P(0)\right)/\beta _{n}(0)$. We choose $\tau '$ small enough such that $\Omega \cap B_{\tau'}\subset \subset \Omega \cup \Gamma$ and $\beta_n\geq \delta_0/2$ on $\Gamma\cap B_{\tau'}$. For $ \varepsilon > 0$, let $\psi_{\varepsilon } (x):=P (x',\varphi(x'))+(A+\eta/2)(x_n-\varphi(x'))+N(x_n-\varphi(x'))^2-\varepsilon |x|^2+\varepsilon \tau '^2/2$. By the definition of $w_{\ast}$ and the continuity of $\psi_{\varepsilon }$, there exists $x^{\ast }\in \bar{\Omega }\cap B_{\tau'}$ such that
\begin{equation}\label{e1.40}
  w(x^{\ast })<w_{\ast}(0)+\frac{1}{4}\varepsilon \tau '^2=\psi_{\varepsilon } (0)-\frac{1}{4}\varepsilon \tau '^2\leq \psi_{\varepsilon } (x^{\ast }).
\end{equation}
By taking $\tau'$ small and $N$ large (as in \Cref{th3.11}), we have
\begin{equation}\label{e2.4}
  F(D^2\psi _{\varepsilon })\geq f~~\mbox{in}~~\Omega \cap B_{\tau '}.
\end{equation}
Similarly, by the continuity of $DP$, $\beta$, $\gamma$ and $g$,
\begin{equation}\label{e2.5}
  \begin{aligned}
    \beta \cdot D\psi _{\varepsilon }+\gamma \psi_{\varepsilon }\geq g~~\mbox{on}~~\Gamma\cap B_{\tau'},
  \end{aligned}
\end{equation}
for $\tau'$ and $\varepsilon$ small enough.

Define
\begin{equation*}
  \bar{w}=\left\{
  \begin{aligned}
    \max & (w,\psi_{\varepsilon } )  &&\mbox{if}~~x\in \bar{\Omega }\cap B_{\tau '};\\
    w &  &&\mbox{if}~~x\notin \bar{\Omega }\cap B_{\tau '}.
  \end{aligned}
  \right.
\end{equation*}
From\cref{e2.4} and\cref{e2.5}, we have that $\psi _{\varepsilon }$ is a subsolution in $\bar{\Omega }\cap B_{\tau '}$. Note that
\begin{equation*}
  w\geq w_{\ast }\geq \psi_{\varepsilon }~~\mbox{on}~~\partial B_{\tau'}\cap \Omega .
\end{equation*}
Thus, $\bar{w}$ is a subsolution of\cref{existence} and hence $\bar{w}\in \textbf{V}$. This implies $\bar{w}\leq w$ which contradicts with\cref{e1.40}.

From above arguments, we conclude that $w_{\ast}$ is a supersolution of\cref{existence}.~\\

\emph{Step 3}. We construct the barriers on the Dirichlet boundary.

Recall\cref{existence_4}. Without loss of generality, we assume that $\xi=e_n$. Then $\beta_n\geq \delta_1$ on $\Gamma$. Given $x_1\in\partial \Omega \backslash \Gamma$, if $x_1 \notin \bar{\Gamma}$, then there exists $r>0$ such that $B_{r}(x_1)\cap \Gamma =\emptyset $. Since $\Omega $ satisfies the exterior cone condition at $x_1$, there exists a function $v_1$ such that $v_1(x_1)=0$, $v_1>0$ in $\Omega\backslash \{x_1\}$ and $F(D^2v_1)\leq -1$ in $\Omega\cap B_{r}(x_1)$.

On the other hand, let $v_2=-K_1x_n^2-K_2x_n+K_3$. By choosing proper constants $K_1$, $K_2$ and $K_3$, we have
\begin{equation}\label{existence:3}
\left\{ \begin{aligned}
&F(D^2v_2)\leq -\|f\|_{L^{\infty }(\Omega ) } &&\mbox{in}~~\Omega;\\
&\beta\cdot Dv_2+\gamma v_2 \leq -\|g\|_{L^{\infty }}-\|\gamma \|_{L^{\infty }}\|\varphi  \|_{L^{\infty }} &&\mbox{on}~~\Gamma;\\
&v_2\geq 1 &&\mbox{in}~~\Omega .
\end{aligned}\right.
\end{equation}
Take $K$ large enough such that $Kv_1>v_2$ on $\Gamma$ and let $v=\inf \{Kv_1,v_2\}$. Then $v$ is a supersolution of\cref{existence} with $v(x_1)=0$ and $v>0$ on $\Omega\backslash {x_1}$.

Since $\varphi $ is continuous at $x_1$, for any $\varepsilon>0$, there exists a constant $k$ large enough such that
\begin{equation}\label{e2.7}
 \begin{aligned}
 w_1:=& \varphi (x_1)-\varepsilon -kv\leq \varphi (x)\\
   \leq& \varphi (x_1)+\varepsilon +kv:=w_2 ~~\mbox{on}~~\partial \Omega \backslash \Gamma.
 \end{aligned}
\end{equation}
Recall\cref{existence:3} and it is easy to verify that $w_1$ is a subsolution and $w_2$ is a supersolution. Hence, $w_1\in \textbf{V}$ and then we have
\begin{equation*}
  w_1\leq w_{\ast}\leq w~~\mbox{on}~~\Omega .
\end{equation*}
By the A-B-P maximum principle and the uniqueness result \Cref{th3.11},
 \begin{equation*}
 w = w^{\ast}\leq w_2~~\mbox{on}~~\Omega .
\end{equation*}
Thus
\begin{equation*}
  w_1\leq w_{\ast}\leq w= w^{\ast}\leq w_2~~\mbox{on}~~\Omega .
\end{equation*}
Since $\varepsilon $ is arbitrary, it follows that $w_{\ast}(x_1)= w(x_1)= w^{\ast}(x_1)=\varphi (x_1)$.

Now, we consider the case $x_1 \in \bar{\Gamma}$. Since $\Omega$ satisfies the exterior sphere condition at $x_1$, let $B_{r_1}(y)$ ($r_1<1$) be the ball such that $B_{r_1}(y)\cap \bar{\Omega}=\{x_1\}$. Let $v_3(x)=r_{1}^{-p}-|x-y|^{-p}$ and $v=\inf\{Kv_3,v_2\}$ ($v_2$ is as above). Similar to above arguments, we conclude that $w_{\ast}(x_1)= w(x_1)= w^{\ast}(x_1)=\varphi (x_1)$.

In consequence, we have
\begin{equation*}
  w_{\ast}= w= w^{\ast}~~\mbox{on}~~\partial \Omega \backslash \Gamma.
\end{equation*}
From the A-B-P maximum principle, we conclude that
\begin{equation*}
  w_{\ast}= w= w^{\ast}~~\mbox{on}~~\bar{\Omega }
\end{equation*}
and hence, $w\in C(\bar{\Omega })$ is a viscosity solution. ~\qed~\\

As a special case of \Cref{th.existence}, the existence of viscosity solutions in a spherical cap is presented in the following. This existence will be used to construct auxiliary functions in later sections.
\begin{corollary}\label{th2.7}
There exists a unique viscosity solution $u\in C(\bar{B}^{+ }_{1,h_0})$ of
\begin{equation}\label{e2.6}
    \left\{
    \begin{aligned}
      &F(D^2u)=0 &&\mbox{in}~~B^{+ }_{1,h_0}; \\
      &\beta \cdot Du =0 &&\mbox{on}~~T_1;\\
      &u=\varphi &&\mbox{on}~~\partial B^{+ }_{1,h_0}\backslash T_1,
    \end{aligned}
    \right.
\end{equation}
where $\beta \in C^{2}(\bar{T}_1)$, $\varphi\in C(\partial B^+_{1,h_0}\backslash T_1)$ and $h_0 >0$ is small enough such that
  \begin{equation}\label{e1.34}
    \beta(x) \cdot\vec{n}(y)<0~~\forall x\in T_1 ~\mbox{and}~ y\in \partial B^+_{1,h_0}\backslash T_1.
  \end{equation}

Furthermore, if $\varphi\in C^{\alpha}(\partial B^+_{1,h_0}\backslash T_1)$, then $u\in C^{\hat{\alpha}}(\bar{B}^+_{1,h_0})$ and
\begin{equation}\label{e2.3}
  \|u\|_{C^{\hat{\alpha}}(\bar{B}^+_{1,h_0})}\leq C\|\varphi \|_{C^{\alpha}(\partial B^+_{1,h_0}\backslash T_1)},
\end{equation}
where $\hat{\alpha}= \alpha/2$ and $C$ depends only on $n$, $\lambda$, $\Lambda$, $\delta_0$ and $h_0$.
\end{corollary}

\begin{remark}\label{r2.8}
From now on, unless stated otherwise, $h_0$ is always chosen (depending only on $\delta_0$) such that\cref{e1.34} holds.
\end{remark}
\proof  We only need to prove the H\"{o}lder estimate. Given $x_1\in\partial  B_{1,h_0}^+\backslash T_1$, since $\varphi\in C^{\alpha}(\partial B^+_{1,h_0}\backslash T_1)$, then (recall\cref{e1.34}) $K\|\varphi \|_{C^{\alpha}(\partial B^+_{1,h_0}\backslash T_1)}(\vec{n}(x_1)\cdot (x-x_1))^{\frac{\alpha }{2}}+\varphi (x_1)$ and $-K\|\varphi \|_{C^{\alpha}(\partial B^+_{1,h_0}\backslash T_1)}(\vec{n}(x_1)\cdot (x-x_1))^{\frac{\alpha }{2}}+\varphi (x_1)$ are supersolution and subsolution respectively for large $K$ depending only on $n$, $\lambda$, $\Lambda$ and $\delta_0$.

From the A-B-P maximum principle, we obtain
\begin{equation*}
  |u(x)-\varphi(x_1)|\leq K\|\varphi \|_{C^{\alpha}(\partial B^+_{1,h_0}\backslash T_1)}|\vec{n}(x_1)\cdot (x-x_1)|^{\frac{\alpha }{2}}\leq K\|\varphi \|_{C^{\alpha}(\partial B^+_{1,h_0}\backslash T_1)}|x-x_1|^{\frac{\alpha }{2}},
\end{equation*}
which implies the H\"{o}lder continuous up to the Dirichlet boundary. Combining with the H\"{o}lder continuous up to the oblique boundary (see \Cref{th5.1}), we obtain the full result\cref{e2.3}.~\qed~\\

\section{ Regularity for the model problem}\label{S4}
In this section, we derive $C^{1,\alpha} $ and $C^{2,\alpha} $ regularity for the model problem, i.e., the homogenous equations in a spherical cap. In addition to the importance of themselves, these regularity will be used to approximate the solution of\cref{e1.1} in different scales and attack the regularity of the solution. Precisely, the model problem is
\begin{equation}\label{model}
    \left\{
    \begin{aligned}
      &F(D^2u)=0 &&\mbox{in}~~ B^+_{1,h_0}; \\
      &\beta \cdot Du = 0 &&\mbox{on}~~T_1,
    \end{aligned}
    \right.
\end{equation}
where $\beta$ is a constant vector with $|\beta |=1$.

From the uniqueness result obtained in last section, the following $C^{1,\alpha}$ regularity for the model problem is derived:
\begin{theorem}\label{th3.5}
Suppose that $u$ is a viscosity solution of\cref{model}. Then $u\in C^{1,\alpha_1 }(\bar{B}^{+}_{1/2,h_0/2})$ and
\begin{equation}\label{e1.5}
  \|u\|_{C^{1,\alpha_1 }(\bar{B}^+_{1/2,h_0/2})}\leq C\left(\|u\|_{L^{\infty }(B^+_{1,h_0})}+ |F(0)|\right),
\end{equation}
where $0<\alpha_1<1 $ depends only on $n$, $\lambda$, $\Lambda$ and $\delta_0$, and $C$ depends also on $h_0$.
\end{theorem}
\proof  Let $v(x)=(u(x+te_i)-u(x))/t^{\alpha_0}$ where $0<t<1/4$ and $i<n$. From \Cref{th3.11}, we have
\begin{equation*}
    \left\{
    \begin{aligned}
     &v\in S (\lambda/n,\Lambda) &&\mbox{in}~~B^+_{1/2,h_0/2}; \\
     &\beta \cdot Dv= 0 &&\mbox{on}~~T_{1/2}.
    \end{aligned}
    \right.
\end{equation*}
By \Cref{th5.1}, $v\in C^{\alpha_0}(\bar{B}^+_{1/4,h_0/4})$ and
\begin{equation*}
\begin{aligned}
  \|v\|_{C^{\alpha _0}(\bar{B}^+_{1/4,h_0/4}) }&\leq  C\|v\|_{L^{\infty }(B^+_{1/2,h_0/2}) }\\
  &\leq C\|u\|_{C^{\alpha _0}(\bar{B}^+_{3/4,3h_0/4}) } \\
  &\leq C(\|u\|_{L^{\infty }(B^+_{1,h_0})}+|F(0)|) (\mbox{again by \Cref{th5.1}}).
\end{aligned}
\end{equation*}
Hence, we obtain that $u$ is $ C^{2\alpha_0}$ along the horizontal directions. Then let $v(x)=(u(x+te_i)-u(x))/t^{2\alpha_0}$ and repeat above procedure, we obtain the $C^{1,\alpha}$ estimate of $u$ on $T_1$ (see \cite[Section 5.3]{C-C} and \cite[Theorem 6.1]{M-S}). Then from the boundary $C^{1,\alpha}$ regularity for fully nonlinear elliptic equations with the Dirichlet boundary conditions,\cref{e1.5} follows.\qed~\\

Next, we intend to derive a boundary $C^{2,\alpha}$ estimate by a similar argument in \cite{M-S}. We introduce the following boundary $C^{1,\alpha}$ estimate for Dirichlet problem which is first proved essentially by Krylov \cite{Kr} and simplified by Caffarelli (see \cite[Theorem 9.32]{G-T}, \cite[Theorem 4.8]{Ka} and \cite[Lemma 7.1]{M-S}).

\begin{proposition}\label{pr3.6}
Suppose that
\begin{equation*}
    \left\{
    \begin{aligned}
      &u\in S(\lambda ,\Lambda ) &&\mbox{in}~~ B^+_{1,h_0}; \\
      &u = 0 &&\mbox{on}~~T_1.
    \end{aligned}
    \right.
\end{equation*}
Then there exists a $C^{\alpha }$ function $A:T_1\rightarrow R$ such that
  \begin{equation*}
    |u(x)-A(x')x_n|\leq C\|u\|_{L^{\infty }(B^+_{1,h_0})}|x_n|^{1+\alpha }~~\forall x\in\bar{B}^+_{1/2,h_0/2}
  \end{equation*}
and
\begin{equation*}
\|A\|_{C^{\alpha }(\bar{T}_{1/2})}\leq C \|u\|_{L^{\infty }(B^+_{1,h_0})},
\end{equation*}
where $\alpha $ depends only on $n$, $\lambda$ and $\Lambda $, and $C$ depends also on $h_0$.
\end{proposition}

To apply above lemma to oblique derivative problems, we rewrite it as follows.
\begin{lemma}\label{le3.7}
Let $u$ be as in \Cref{pr3.6}. Then there exists a $C^{\alpha }$ function $\bar{A}:T _1\rightarrow R$ such that
  \begin{equation*}
    |u(x)-\bar{A}(x_{\beta }')|x-(x_{\beta }',0)||\leq C\|u\|_{L^{\infty }(B^+_{1,h_0})}|x-(x_{\beta }',0)|^{1+\alpha }~~\forall x\in\bar{B}^{+ }_{1/2,h_0/2},
  \end{equation*}
and
\begin{equation*}
  \|\bar{A}\|_{C^{\alpha }(\bar{T}_{1/2})}\leq C \|u\|_{L^{\infty }(B^+_{1,h_0})},
\end{equation*}
where $\alpha$ and $C$ are as in \Cref{pr3.6}.
\end{lemma}
\begin{remark}\label{}
  Note that for every $x\in B^+_{1,h_0}$, $(x'_{\beta },0)\in T_1$.
\end{remark}
\proof  Let $\bar{A}(x_{\beta }')=A(x_{\beta }')\beta _n$. Then by \Cref{pr3.6},
\begin{equation*}
\begin{aligned}
   |u(x)-\bar{A}(x_{\beta }')|x-(x_{\beta }',0)||=&|u(x)-A(x_{\beta }')\beta _n\cdot \frac{x_n}{\beta _n}| \\
  = & |u(x)-A(x')x_n+A(x')x_n-A(x_{\beta }')x_n|\\
  \leq & C\|u\|_{L^{\infty }(B^{+ }_{1,h_0})}|x_n|^{1+\alpha }+C\|A\|_{C^{\alpha }(\bar{T}_{1/2})}|x'-x_{\beta }'|^{\alpha }x_n\\
  \leq & C\|u\|_{L^{\infty }(B^{+ }_{1,h_0})}|x-(x_{\beta }',0)|^{1+\alpha }~~\forall x\in\bar{B}^{+ }_{1/2,h_0/2}.
\end{aligned}
\end{equation*}\qed~\\

\begin{lemma}\label{co3.8}
Suppose that $u$ is a viscosity solution of\cref{model}. Then there exists a $C^{\alpha }$ function $\bar{A}:T _1\rightarrow R$ such that for any $x\in\bar{B}^+_{1/2,h_0/2}$,
\begin{equation*}
    |u(x)-u(x_{\beta }',0)-\frac{1}{2}\bar{A}(x_{\beta }')|x-(x_{\beta }',0)|^2|\leq C(\|u\|_{L^{\infty }(B^+_{1,h_0})}+|F(0)|)|x-(x_{\beta }',0)|^{2+\alpha }
\end{equation*}
and
\begin{equation*}
  \|\bar{A}\|_{C^{\alpha }(\bar{T}_{1/2})}\leq C (\|u\|_{L^{\infty }(B^+_{1,h_0})}+|F(0)|),
\end{equation*}
where $\alpha$ is as in \Cref{pr3.6} and $C$ depends only on $n$, $\lambda$, $\Lambda$, $\delta_0$ and $h_0$ .
\end{lemma}
\proof  By \Cref{th3.5}, $u\in C^{1,\alpha_1 }(\bar{B}^{+ }_{7/8,7h_0/8})$. Hence, $u_{\beta}:=\beta \cdot Du\in C^{\alpha }(\bar{B}^{+ }_{7/8,7h_0/8})$. From the closedness of viscosity solutions (\Cref{le2.3}), $u_{\beta }\in S(\lambda /n,\Lambda )$ in $B^{+ }_{7/8,7h_0/8}$ and $u_{\beta }=0$ on $T_{7/8}$. Then by \Cref{le3.7}, there exists $\bar{A}$ such that
\begin{equation*}
  \begin{aligned}
     u(x)-u(x_{\beta }',0)= & \int_{0}^{|x-(x_{\beta }',0)|} u_{\beta }\left((x_{\beta }',0)+ t\beta\right)dt\\
     \leq & \int_{0}^{|x-(x_{\beta }',0)|} \left(\bar{A}(x_{\beta }')t+ C\|u_{\beta }\|_{L^{\infty }({B^{+ }_{7/8,7h_0/8})}}t^{1+ \alpha }\right) dt\\
     =& \frac{\bar{A}(x_{\beta }')}{2} |x-(x_{\beta }',0)|^{2}+C(\|u\|_{L^{\infty }(B^+_{1,h_0})}+|F(0)|)|x-(x_{\beta }',0)|^{2+\alpha }.
   \end{aligned}
\end{equation*}
Hence,
\begin{equation*}\label{e1.6}
 u(x)-u(x_{\beta }',0)-\frac{1}{2}\bar{A}(x_{\beta }')|x-(x_{\beta }',0)|^2\leq C(\|u\|_{L^{\infty }(B^+_{1,h_0})}+|F(0)|)|x-(x_{\beta }',0)|^{2+\alpha }.
\end{equation*}
The proof of the other direction of the inequality is similar and we omit it. Finally, note that
\begin{equation*}
  \|\bar{A}\|_{C^{\alpha }(\bar{T}_{1/2})}\leq C\|u_{\beta }\|_{L^{\infty }(B^+_{7/8,7h_0/8})}\leq  C(\|u\|_{L^{\infty }(B^+_{1,h_0})}+|F(0)|)
\end{equation*}
and we complete the proof. \qed~\\

Next, we follow the idea of \cite{M-S} to show that the solution satisfies an equation on the flat boundary.
\begin{lemma}\label{le3.9}
Suppose that $u$ is a viscosity solution of\cref{model}. Define $v:T_1\rightarrow R $ by $v(x')=u(x',0)$. Then
\begin{equation}\label{e1.7}
  G(D^2v,x')=0~~\mbox{in}~~T_1
\end{equation}
where
\begin{equation*}
  G(M,x'):=F\left(\left(
                   \begin{array}{cc}
                     M & -M\beta '/\beta _{n} \\
                     - \beta '^{T}M/\beta _{n}& ~~ \frac{1}{\beta _{n}^{2}} \left(\bar{A}(x')+\beta '^{T}M\beta '\right)\\
                   \end{array}
                 \right)\right)
\end{equation*}
for any $x'\in T_1$ and $M\in S^{n-1}$.
\end{lemma}
\proof  Let $\varphi \in C^2(T_1)$ touch $v$ by below at $x_0\in T_1$. Without loss of generality, we assume that $x_0=0$. From \Cref{co3.8}, we have
\begin{equation}\label{e1.8}
u(x)\geq u(x_{\beta }',0)+\frac{1}{2}\bar{A}(x_{\beta }')|x-(x_{\beta }',0)|^2- C(\|u\|_{L^{\infty }(B^{+ }_{1,h_0})}+|F(0)|)|x-(x_{\beta }',0)|^{2+\alpha }.
\end{equation}
Then
\begin{equation*}
  \begin{aligned}
    u(x)\geq & u(x_{\beta }',0)+\frac{1}{2}\bar{A}(x_{\beta }')|x-(x_{\beta }',0)|^2- C|x-(x_{\beta }',0)|^{2+\alpha } \\
    \geq & \varphi (x_{\beta }')+\frac{1}{2}\bar{A}(x_{\beta }')|x-(x_{\beta }',0)|^2- C|x-(x_{\beta }',0)|^{2+\alpha }\\
    =& \varphi (x'-\beta'\frac{x_n}{\beta_n})+ \frac{1}{2}\bar{A}(x_{\beta }')\frac{x_n^2}{\beta_n^2}-Cx_n^{2+\alpha }\\
    \geq & \varphi (x'-\beta'\frac{x_n}{\beta_n})+ \frac{1}{2}\bar{A}(0)\frac{x_n^2}{\beta_n^2}-C|x_{\beta }'|^{\alpha }x_n^{2}-Cx_n^{2+\alpha }.
  \end{aligned}
\end{equation*}
Let $\varphi _{\varepsilon }(x',x_n)=\varphi (x'-(x_n/\beta_n)\beta')+\bar{A}(0)x_n^2/(2\beta_n^2)-2\varepsilon x_n^{2}$. Then
\begin{equation}\label{e3.3}
  u\geq \varphi _{\varepsilon }+\varepsilon x_n^2~~\mbox{in}~~B^+_r~~\mbox{for}~~r~~\mbox{small enough and } u(0)=\varphi _{\varepsilon }(0)=\varphi (0).
\end{equation}
Let $\bar{\varphi}_{\varepsilon } (x)=\varphi _{\varepsilon }(x+t\beta )$ and $y=x'+t\beta'-(x_n+t\beta_n)\beta'/\beta_n$. Then
\begin{equation*}
\begin{aligned}
D\bar{\varphi}_{\varepsilon } (x)= \left(D\varphi \left(y\right),\frac{1}{\beta_n^{2}}\bar{A}(0)(x_n+t\beta_n)-
4\varepsilon(x_n+t\beta_n)-\beta'\cdot D\varphi \left(y\right)\frac{1}{\beta_n}\right).
\end{aligned}
\end{equation*}
Hence,
\begin{equation*}
  \beta \cdot D\bar{\varphi}_{\varepsilon } (x',0)=\left(\bar{A}(0)-4\varepsilon \beta _{n}^{2}\right)t.
\end{equation*}
If $\bar{A}(0)\leq 0$, then $\beta \cdot D\bar{\varphi}_{\varepsilon } (x',0)>0$ on $T_r$. We take $t<0$ small enough. Hence, $\bar{\varphi}_{\varepsilon }+c_0$ will touch $u$ by below at some $x_0\in B^{+}_{r/2}$ for a proper $c_0$. Since $F(D^2u)\leq 0$, we have
\begin{equation*}
F\left(\left(
                   \begin{array}{cc}
                     D^2\varphi  & -D^2\varphi \beta '/\beta _{n} \\
                     - \beta '^{T}D^2\varphi /\beta _{n}& ~~ \frac{1}{\beta _{n}^{2}} \left(\bar{A}(0)+\beta '^{T}D^2\varphi \beta '\right)-4\varepsilon\\
                   \end{array}
                \right) \right)\leq 0.
\end{equation*}
Let $\varepsilon \rightarrow 0$, by the continuity of $D^2\varphi$, we have
\begin{equation}\label{e1.10}
F\left(\left(
                   \begin{array}{cc}
                     D^2\varphi (0) & -D^2\varphi (0) \beta '/\beta _{n} \\
                     - \beta '^{T}D^2\varphi (0) /\beta _{n}& ~~ \frac{1}{\beta _{n}^{2}} \left(\bar{A}(0)+\beta '^{T}D^2\varphi (0) \beta '\right)\\
                   \end{array}
                \right) \right)\leq 0.
\end{equation}
If $\bar{A}(0)> 0$, we take $\varepsilon < \bar{A}(0)/(8\beta _{n}^{2})$. Then, for $r$ small enough, we also have\cref{e3.3}. By taking $t>0$ small, $\bar{\varphi}_{\varepsilon }+c_0$ will touch $u$ by below at some $x_1\in B^{+}_{r/2}$ for a proper $c_0$. Similarly,
\begin{equation*}
F\left(\left(
                   \begin{array}{cc}
                     D^2\varphi  & -D^2\varphi \beta '/\beta _{n} \\
                     - \beta '^{T}D^2\varphi /\beta _{n}& ~~ \frac{1}{\beta _{n}^{2}} \left(\bar{A}(0)+\beta
                      '^{T}D^2\varphi \beta '\right)-4\varepsilon\\
                   \end{array}
                 \right)\right)\leq 0.
\end{equation*}
Let $\varepsilon \rightarrow 0$, we also obtain\cref{e1.10}. Therefore, $v$ is a supersolution of $G(D^2v,x')=0$. The verification of subsolution is similar and we omit it. \qed~\\

Now we prove the boundary $C^{2,\alpha }$ estimate.
\begin{theorem}\label{th3.10}
Suppose that $F$ is convex and $u$ is a viscosity solution of\cref{model}. Then $u\in C^{2,\alpha_2 } (\bar{B}^+_{1/2,h_0/2})$ and
\begin{equation}\label{e1.11}
   \|u\|_{C^{2,\alpha_2 }(\bar{B}^{+ }_{1/2,h_0/2})}\leq C\left(\|u\|_{L^{\infty }(B^+_{1,h_0})}+ |F(0)|\right),
\end{equation}
where $0<\alpha _2<1$ depends only on $n$, $\lambda$, $\Lambda$ and $\delta_0$, and $C$ depends also on $h_0$.
\end{theorem}
\proof  By \Cref{le3.9}, $v(x')=u(x',0)$ satisfies $G(D^2v,x')=0$ in $T_1$. Note that $G$ is uniformly elliptic with $\lambda$ and $\tilde{\Lambda}$ where $\tilde{\Lambda}$ depends only on $n $, $\Lambda$ and $\delta_0 $. Furthermore, $G$ is convex. Indeed,
\begin{equation*}
  \begin{aligned}
    G&\left (tM+(1-t)N,x'\right )  \\
 &=F\left(t\left(
                   \begin{array}{cc}
                     M & -M\beta '/\beta _{n} \\
                     - \beta '^{T}M/\beta _{n}& \frac{1}{\beta _{n}^{2}} \left(\bar{A}(x')+\beta '^{T}M\beta '\right)\\
                   \end{array}
                 \right)\right.\\
&~~~~+\left.(1-t) \left(
                   \begin{array}{cc}
                     N & -N\beta '/\beta _{n} \\
                     - \beta '^{T}N/\beta _{n}& \frac{1}{\beta _{n}^{2}} \left(\bar{A}(x')+\beta '^{T}N\beta '\right)\\
                   \end{array}
                 \right)
                 \right)\\
&\leq tG(M,x')+(1-t)G(N,x').
  \end{aligned}
\end{equation*}
Note that $G$ is H\"{o}lder continuous in $x'$. By the classical interior $C^{2,\alpha}$ estimates (Evans-Krylov estimates) for convex operators, there exists a constant $0<\alpha<1$ depending only on $n$, $\lambda$, $\Lambda$ and $\delta_0$ such that $v\in C^{2,\alpha}(T_{3/4})$ and
\begin{equation*}
   \|v\|_{C^{2,\alpha }(\bar{T}_{3/4})}\leq C\left(\|v\|_{L^{\infty }(T_1)}+ |F(0)|\right).
\end{equation*}
Note that $u=v$ on $T_1$. By the boundary $C^{2,\alpha }$ estimates for Dirichlet problems, there exists a constant $0<\alpha_2<1$ depending only on $n$, $\lambda$, $\Lambda$ and $\delta_0$ such that $u\in C^{2,\alpha_2 }(\bar{B}^+_{1/2,h_0/2})$ and\cref{e1.11} holds. \qed~\\
\begin{remark}\label{re4.4zk}
  From now on, $\alpha_0$, $\alpha_1$ and $\alpha_2$ always denote the constants originated from \Cref{th5.1}, \Cref{th3.5} and \Cref{th3.10}, and depending only on $n$, $\lambda$, $\Lambda$ and $\delta_0$.
\end{remark}

\section{\texorpdfstring{$C^{1,\alpha}$}{C1,alpha} regularity}\label{S5}
In the following sections, we use perturbation method to deduce the $C^{1,\alpha}$, $C^{2,\alpha}$ and higher regularity for oblique derivative problem\cref{e1.1}. The main idea of perturbation method is applying the solutions of the model problem\cref{model} to approximate the solution of\cref{e1.1}. The existence of solutions is assured by \Cref{th2.7} and the sufficient regularity has been obtained in last section. In approximating the solution of\cref{e1.1}, $f(x)$, $\beta(x)$ and $\Gamma$ are regarded as the perturbation of $0$, $\beta(0)$ and $T_1$ etc.. The A-B-P maximum principle is the main tool to measure the difference between the solution of\cref{e1.1} and the solutions of the model problem. In this section, we deduce the $C^{1,\alpha }$ regularity for\cref{e1.1}.
\begin{lemma}\label{thh4.1}
Let $u$ be a viscosity solution of
  \begin{equation*}
    \left\{
    \begin{aligned}
      &F(D^2u)=f &&\mbox{in}~~ \Omega; \\
      &\beta  \cdot Du = g &&\mbox{on}~~\Gamma,
    \end{aligned}
    \right.
  \end{equation*}
$x_0\in \Gamma$ such that $\mathrm{dist}(x_0,\partial\Omega\backslash \Gamma)>1$ and $0<\alpha <\alpha _1$.

Suppose that $\beta$ and $g$ are $C^{\alpha}$ at $x_0\in \Gamma$, and $f$ satisfies
\begin{equation}\label{e1.12}
\left (\frac{1}{|B_r(x_0)\cap\Omega|}\int_{B_r(x_0)\cap\Omega}|f|^n\right )^{\frac{1}{n}}\leq C_f r^{\alpha - 1}~~\forall r>0.
\end{equation}
Then $u$ is $C^{1,\alpha}$ at $x_0$, i.e., there exists an affine function $l$ such that
\begin{equation}\label{e1.15}
   \|u-l\|_{L^{\infty }(\bar{\Omega }\cap B_r(x_0))}\leq C_0r^{1+\alpha} ~~~~\forall 0<r<r_1,
\end{equation}
\begin{equation}\label{e4.49}
  |Dl|\leq C_0,
\end{equation}
\begin{equation}\label{e4.50}
  C_0\leq C\left(\|u\|_{L^{\infty }(\Omega )}+C_f+\|g\|_{C^{\alpha }(x_0)}+|F(0)|\right)
\end{equation}
and
\begin{equation}\label{e4.51}
  r_1=\check{C}^{-1},
\end{equation}
where $C$ depends only on $n$, $\lambda$, $\Lambda$, $\delta_0$ and $\alpha$, and $\check{C}$ depends also on $\|\beta \|_{C^{\alpha }(x_0)}$ and the $C^{1}$ modulus of $\Gamma$ at $x_0$.
\end{lemma}
\begin{remark}\label{re5.1.zk}
As in \Cref{th5.1}, the condition $\mathrm{dist}(x_0,\partial\Omega\backslash \Gamma)>1$ is not an essential assumption and ``1'' can be replaced by any positive constant. Then, we obtain the scaling version of\crefrange{e1.15}{e4.51}.
\end{remark}
\proof
We make some normalization first. We assume that $F(0)=0$. Otherwise, by the uniform ellipticity, there exists $t\in R$ such that $F(tI)=0$ and $|t|\leq  |F(0)|/\lambda$. Then $\tilde{u}=u-t|x|^2/2$ satisfies
 \begin{equation*}
    \left\{
    \begin{aligned}
      &F(D^2\tilde{u}+tI)=f &&\mbox{in}~~ \Omega; \\
      &\beta \cdot D\tilde{u} = g-\frac{t}{2}\beta  \cdot x &&\mbox{on}~~\Gamma.
    \end{aligned}
    \right.
  \end{equation*}
Hence, the estimate of $u$ follows from that of $\tilde{u}$ . Next, by choosing a proper coordinate, we assume that $x_0$ is the origin,
\begin{equation*}\label{e1.134}
\Gamma \cap B_{1}=\{(x',x_n)\in B_1\big|x_n=\varphi (x')\}
\end{equation*}
and $\varphi $ satisfies that
\begin{equation*}\label{e1.144}
\varphi(0)=0,~~D\varphi(0)=0~~\mbox{and}~~ |\varphi (x')|\leq \nu|x'|,
\end{equation*}
where $\nu$ is chosen small later.

Denote $\beta(0)$ by $\beta^0$. Since $\beta$ is $C^{\alpha }$ at $0$ and $\beta^0_n=\beta(0)\cdot \vec{n}(0)\geq \delta_0$, $\beta_n\geq \delta_0 /2$ on $\Gamma \cap B_{r_0}$ for $r_0$ small enough (depending only on $\delta_0$ and $\|\beta \|_{C^{\alpha }(x_0)}$). Without loss of generality, we assume that $r_0=1$. By scaling, we also assume that $[\beta ]_{C^{\alpha } (0)}\leq 1$. Finally, we assume that $g(0)=0$. Otherwise, we may consider $\tilde{u}=u-g(0)x_n/\beta^0_n$.

Let $M:=\|u\|_{L^{\infty }(\Omega )}+C_f+\|g\|_{C^{\alpha }(0)}$ and $\Omega _{r}:=\Omega \cap B_{r,h_0 r}$ where $h_0$ is chosen as in \Cref{th2.7} (depending only on $\delta_0$) such that\cref{e1.34} holds with $\beta^0$. To prove that $u$ is $C^{1,\alpha}$ at $0$, we only need to prove the following:

There exist constants $0< \tau <1$, $0< \eta < 1$, $\bar{C}$ and $\hat{C}$ depending only on $n$, $\lambda$, $\Lambda$, $\delta_0$ and $\alpha$, and a sequence affine functions $l_k(x)=b_kx+c_k$ ($k\geq -1$) such that for all $k\geq 0$

\begin{equation}\label{e1.16}
\|u-l_k\|_{L^{\infty }(\Omega _{\tau \eta^{k}})}\leq \hat{C} M \eta ^{k(1+\alpha )},
\end{equation}

\begin{equation}\label{e1.17}
\eta ^{k}|b_k-b_{k-1}|+ |c_k-c_{k-1}|\leq \bar{C}\hat{C}M\eta ^{k(1+\alpha )}
\end{equation}
and
\begin{equation}\label{e1.18}
  \beta^0\cdot b_k=g(0)=0.
\end{equation}

We prove above by induction. For $k=0$, by setting $l_0=l_{-1}=0$, the conclusion holds clearly. Suppose that the conclusion holds for $k=k_0$. We need to prove that the conclusion holds for $k=k_0+1$. In the rest of the proof, $C$, $C_1$, $C_2$ etc. denote positive constants depending only on $n$, $\lambda$, $\Lambda$, $\delta_0$ and $\alpha$.

Let $r:=\frac{1}{2}\tau \eta ^{k_{0}}$, $\tilde{B}^{+}_{r,h_0r}:=B^{+}_{r,h_0r}-\nu r e_n $, $\tilde{T}_r:=T_r-\nu r e_n$ and $\tilde{\Omega }_{r}:=\Omega \cap \tilde{B}^{+}_{r,h_0r}$. Assume that
\begin{equation}\label{e1.19}
  4\nu \leq h_0.
\end{equation}
Then $\Omega _{r/2}\subset \tilde{\Omega }_{r}$.

Note that $u-l_{k_0}$ satisfies
 \begin{equation*}
    \left\{
    \begin{aligned}
      & F(D^2(u-l_{k_0}))=f &&\mbox{in}~~\Omega _{2r};\\
      &\beta \cdot D(u-l_{k_0}) = g &&\mbox{on}~~\Gamma\cap \partial \Omega _{2r}.
    \end{aligned}
    \right.
  \end{equation*}
By the H\"{o}lder estimate (\Cref{th5.1}),
\begin{equation}\label{e.21}
\begin{aligned}
   &\|u-l_{k_0}\|_{L^{\infty }(\Omega _{r})}+r^{\alpha _0}[u-l_{k_0}]_{C^{\alpha _0}(\bar{\Omega} _{r})}\\
   &\leq C\left (\|u-l_{k_0}\|_{L^{\infty }(\Omega _{2r})}+r\|f\|_{L^n(\Omega_{2r})}+r\|g\|_{L^{\infty }(\Gamma  \cap B_{2r})}\right ) \\
   &\leq C\left (\|u-l_{k_0}\|_{L^{\infty }(\Omega _{2r})}+r^{1+\alpha }C_f+r\|g\|_{L^{\infty }(\Gamma  \cap B_{2r})}\right ).
\end{aligned}
\end{equation}
Extend $u-l_{k_0}$ from $\Omega _{r}$ to the whole $R^n$ such that
\begin{equation}\label{e1.45}
  \|u-l_{k_0}\|_{L^{\infty }(R^n)}+r^{\alpha _0}[u-l_{k_0}]_{C^{\alpha _0}(R^n )}\leq C \left ( \|u-l_{k_0}\|_{L^{\infty }(\Omega _{r})}+r^{\alpha _0}[u-l_{k_0}]_{C^{\alpha _0}(\bar{\Omega} _{r})}\right ).
\end{equation}

By \Cref{th2.7}, there exists a unique solution $v$ of
\begin{equation*}
    \left\{
    \begin{aligned}
      &F(D^2v)=0 &&\mbox{in}~~ \tilde{B}^{+}_{r,h_0r} ; \\
      &\beta^0 \cdot Dv = 0 &&\mbox{on}~~\tilde{T}_r;\\
      &      v= u-l_{k_0} &&\mbox{on}~~\partial \tilde{B}^{+}_{r,h_0r}\backslash\tilde{T}_r.
    \end{aligned}
    \right.
\end{equation*}
Let $w=u-l_{k_0}-v$. Then by \Cref{th3.11}, $w$ satisfies
\begin{equation*}
    \left\{
    \begin{aligned}
      &w\in S(\lambda /n,\Lambda , f) &&\mbox{in}~~ \Omega \cap \tilde{B}^{+}_{r,h_0r}; \\
      &\beta\cdot  Dw = g-\beta\cdot (b_{k_0}+Dv) &&\mbox{on}~~\Gamma  \cap \tilde{B}^{+}_{r,h_0r};\\
            &w= 0 &&\mbox{on}~~\partial \tilde{B}^{+}_{r,h_0r}\cap \bar{\Omega}.
    \end{aligned}
    \right.
\end{equation*}

In the following arguments, we estimate $v$ and $w$ respectively. By the boundary $C^{1,\alpha_1}$ estimates for $v$ (\Cref{th3.5}) and the A-B-P maximum principle (\Cref{th2.4}), there exists an affine function $\tilde{l}(x)=\bar{b}\left(x+\nu re_n\right)+\bar{c}$  such that
\begin{equation*}
\begin{aligned}
  \|v-\tilde{l}\|_{L^{\infty }(\Omega _{2\eta r})}\leq C\frac{(2\eta r)^{1+ \alpha_1 }}{r^{1+ \alpha _{1}}}\|v\|_{L^{\infty }( \tilde{B}^{+}_{r,h_0r})}&\leq C\eta ^{1+ \alpha _1}\|u-l_{k_0}\|_{L^{\infty }(\Omega _{ r})}\\
  & \leq C\eta ^{\alpha _1-\alpha }\cdot \hat{C}M\eta ^{(k_0+1)(1+\alpha)},
\end{aligned}
\end{equation*}
\begin{equation}\label{e.19}
 \eta ^{k_0} |\bar{b}|+ |\bar{c}|\leq \bar{C}\cdot \hat{C}M\eta ^{k_0(1+\alpha)}
\end{equation}
and
\begin{equation}\label{e.23}
  \beta^0\cdot \bar{b}=0.
\end{equation}
Let $\bar{l}(x)=\bar{b}x+\bar{c}$. Then
\begin{equation}\label{e.20}
\begin{aligned}
   \|v-\bar{l}\|_{L^{\infty }(\Omega _{2\eta r})}\leq & \|v-\tilde{l}\|_{L^{\infty }(\Omega _{2\eta r})}+ \|\bar{l}-\tilde{l}\|_{L^{\infty }(\Omega _{2\eta r})}\\
   \leq & C\eta ^{\alpha _1-\alpha }\cdot \hat{C}M\eta ^{(k_0+1)(1+\alpha)}+|\nu \bar{b}_nr|\\
  \leq & \left (C_1\eta ^{\alpha _1-\alpha }+\frac{C_2 \tau}{\eta ^{1+\alpha }} \right )\cdot \hat{C}M\eta ^{(k_0+1)(1+\alpha)}.
\end{aligned}
\end{equation}

Next, we estimate the term $w$. Let $\tilde{B}^{+}_{\mu }:=\{x\in \tilde{B}^{+ }_{r,h_0r}\big|\mathrm{dist}(x,\partial \tilde{B}^{+ }_{r,h_0r}\backslash\tilde{T}_r)\geq \mu r\}$, $\Omega _{r,\mu }:=\Omega \cap \tilde{B}^{+}_{\mu }$, $\Gamma_1:=\partial \tilde{B}^{+}_{\mu }\cap \bar{\Omega} $ and $\Gamma_2:=\tilde{B}^{+}_{\mu }\cap \Gamma $.

By the global H\"{o}lder estimate for $v$ (\Cref{th2.7}) and recalling\cref{e.21} and\cref{e1.45}, there exists $0<\alpha_3\leq \alpha_0/2$ such that
\begin{equation*}\label{e.22}
\begin{aligned}
    \|v\|_{L^{\infty }(\tilde{\Omega} _{r})}+r^{\alpha _3}[v]_{C^{\alpha _3}(\tilde{\Omega} _{r})}\leq & C\left (\|u-l_{k_0}\|_{L^{\infty }(\Omega _{r})}+r^{\alpha _0}[u-l_{k_0}]_{C^{\alpha _0}(\Omega _{r})}\right )\\
    \leq &C\left(\|u-l_{k_0}\|_{L^{\infty }(\Omega _{2r})}+r^{1+\alpha }C_f+r\|g\|_{L^{\infty }(\Gamma  \cap B_{2r})}\right). \\
\end{aligned}
\end{equation*}
For any $ x_0\in \Gamma_1$, there exists $ \bar{x}\in \partial \tilde{B}^{+ }_{r,h_0r}\backslash\tilde{T}_r$ such that $|x_0-\bar{x}|=\mu r$. Then by recalling\cref{e.21},
\begin{equation}\label{e1.28}
\begin{aligned}
    |w(x_0)|=&|u(x_0)-l_{k_0}(x_0)-v(x_0)|\\
    =&|u(x_0)-l_{k_0}(x_0)-v(x_0)-u(\bar{x})+l_{k_0}(\bar{x})+v(\bar{x})| \\
    \leq & |\left (u(x_0)-l_{k_0}(x_0)\right )-\left (u(\bar{x})-l_{k_0}(\bar{x})\right )|+|v(x_0)-v(\bar{x})|\\
    \leq & C\frac{(\mu r)^{\alpha _{3}}}{r^{\alpha _{3}}}\left(\|u-l_{k_0}\|_{L^{\infty }(\Omega _{2r })}+r^{1+\alpha }C_f+r\|g\|_{L^{\infty }(\Gamma  \cap B_{2r})}\right).
\end{aligned}
\end{equation}
Combining the A-B-P maximum principle and\cref{e1.28}, we have
\begin{equation}\label{e1.22}
  \begin{aligned}
    \|w\|&_{L^{\infty }(\Omega _{r,\mu })} \leq \|w\|_{L^{\infty }(\Gamma_1)}+Cr\|g-\beta\cdot  (b_{k_0}+Dv)\|_{L^{\infty }(\Gamma_2)}+r\|f\|_{L^n(\Omega_{2r})}\\
    \leq & C\frac{(\mu r)^{\alpha _{3}}}{r^{\alpha _{3}}}\left(\|u-l_{k_0}\|_{L^{\infty }(\Omega _{2r})}+r^{1+\alpha }C_f+r\|g\|_{L^{\infty }(\Gamma  \cap B_{2r})}\right)\\
    &+Cr\|g\|_{L^{\infty }(\Gamma_2)}+Cr\|\beta \cdot b_{k_0}\|_{L^{\infty }(\Gamma_2)}+Cr\|\beta\cdot  Dv\|_{L^{\infty }(\Gamma_2)}+Cr^{1+\alpha}C_f\\
    \leq& C_3\mu ^{\alpha _3}\hat{C}M\eta ^{k_0(1+\alpha)}+C_4r^{1+\alpha}\|g\|_{C^{\alpha }(\bar{\Gamma}  \cap \bar{B}_{2r})}+Cr\|\beta\cdot  b_{k_0}\|_{L^{\infty }(\Gamma_2)}\\
    &+Cr\|\beta\cdot  Dv\|_{L^{\infty }(\Gamma_2)}+C_5r^{1+\alpha}C_f\\
    \leq & \left (C_3\mu ^{\alpha _3}+\frac{C_4+C_5}{\hat{C}}\right )\hat{C}M\eta ^{k_0(1+\alpha)}+Cr\|\beta \cdot b_{k_0}\|_{L^{\infty }(\Gamma_2)}+Cr\|\beta\cdot  Dv\|_{L^{\infty }(\Gamma_2)}.
  \end{aligned}
\end{equation}

In the following, we estimate $\|\beta \cdot b_{k_0}\|_{L^{\infty }(\Gamma_2)}$ and $\|\beta\cdot  Dv\|_{L^{\infty }(\Gamma_2)}$ respectively. For the first term, recall that $\beta^{0}\cdot b_{k_0}=0$ and then we obtain
\begin{equation}\label{e1.23}
  \begin{aligned}
   \|\beta\cdot  b_{k_0}\|_{L^{\infty }(\Gamma_2)}=&\|(\beta-\beta^{0})\cdot  b_{k_0}\|_{L^{\infty }(\Gamma_2)}\leq \frac{C\hat{C}M}{1-\eta }[\beta ]_{C^{\alpha}(x_0)}r^{\alpha }\leq C_6r^{\alpha }\hat{C}M.
  \end{aligned}
\end{equation}

We assume that
\begin{equation}\label{e1.24}
  4\nu \leq \mu.
\end{equation}
Then, $\forall x_0 \in \Gamma_2$, $\mathrm{dist} (x_0,\tilde{T}_r )< \frac{1}{2}\mathrm{dist} (x_0,\partial \tilde{B}^{+}_{r,h_0r}\backslash \tilde{T}_r )$. Let $x^{\ast }\in \tilde{T}_r$ such that $|x_0-x^{\ast}|=\mathrm{dist} (x_0,\tilde{T}_r )$. By the $C^{1,\alpha _1}$ estimate for $v$ in $B_{\mu r}(x^{\ast})\cap \tilde{B}^{+}_{r,h_0r}$, we have (note that $\beta^{0}\cdot Dv(x^{\ast})=0$)
\begin{equation}\label{e1.25}
  \begin{aligned}
   |\beta(x_0) &\cdot Dv(x_0)|\\
   \leq& |\beta(x_0) -\beta^{0}||Dv(x_0)|+|\beta^{0}||Dv(x_0)-Dv(x^{\ast})|\\
   \leq & \left([\beta ]_{C^{\alpha}(x_0)}\frac{r^{\alpha }}{\mu r}+|\beta^{0}|\frac{(2\nu r)^{\alpha _1}}{(\mu r)^{1+ \alpha _1}}\right)\|v\|_{L^{\infty }(B_{\mu r}(x^{\ast})\cap \tilde{B}^{+}_{r,h_0r})}\\
   \leq & \left([\beta ]_{C^{\alpha}(x_0)}\frac{r^{\alpha - 1}}{\mu}+\frac{(2\nu)^{\alpha _1}|\beta^{0}|}{\mu^{1+ \alpha _1}r}\right)\cdot \hat{C}M\eta ^{k_0(1+\alpha)}\\
   =& \left(C_7\frac{r^{\alpha - 1}}{\mu}+\frac{C_8\nu ^{\alpha _1} }{\mu^{1+ \alpha_1 }r}\right)\cdot \hat{C}M\eta ^{k_0(1+\alpha)}.
  \end{aligned}
\end{equation}
Combining\cref{e1.22},\cref{e1.23} and\cref{e1.25}, we have
\begin{equation}\label{e1.26}
\begin{aligned}
  \|w\|_{L^{\infty }(\Omega _{r,\mu })}\leq \bigg(C_3\mu ^{\alpha _3}+\frac{C_4+C_5}{\hat{C} }+C_6\tau ^{1+ \alpha }+\frac{C_7\tau ^{\alpha }}{\mu} +\frac{C_8\nu ^{\alpha _1} }{\mu^{1+ \alpha_1 }}\bigg)\hat{C}M\eta ^{k_0(1+\alpha)}.
\end{aligned}
\end{equation}

Take $\eta $ small enough such that $ C_1\eta ^{\alpha _1- \alpha }<1/4$. Let $\mu =\tau^{\alpha /2}$ and take $\tau$ small enough such that
\begin{equation*}
\begin{aligned}
   &\frac{C_2\tau}{\eta ^{1+ \alpha }}<\frac{1}{4} ,~\frac{C_3\tau ^{\alpha \alpha _3/2}}{\eta ^{1+ \alpha }}<\frac{1}{12},~\frac{C_6\tau ^{1+\alpha }}{\eta ^{1+\alpha}}<\frac{1}{12} \mbox{ and }\frac{C_7\tau ^{\alpha/2 }}{\eta ^{1+\alpha}}<\frac{1}{12}. \\
\end{aligned}
\end{equation*}
Next, take $\nu$ small enough such that\cref{e1.19} and\cref{e1.24} hold and
\begin{equation*}
\frac{C_8\nu ^{\alpha _1} }{\mu^{1+ \alpha_1 } \eta ^{1+ \alpha }}\leq \frac{1}{12}.
\end{equation*}
Finally, take $\hat{C}$ large enough such that
\begin{equation*}
  \frac{C_4+C_5}{\hat{C} \eta ^{1+\alpha}}\leq \frac{1}{12}.
\end{equation*}
Therefore, combining\cref{e.20} and\cref{e1.26}, we have
\begin{equation*}
\begin{aligned}
  \|u-l_{k_0}-\bar{l}\|_{L^{\infty }(\Omega _{\tau \eta ^{(k_0+1)}})}&=  \|u-l_{k_0}-v+v-\bar{l}\|_{L^{\infty }(\Omega _{\tau \eta ^{(k_0+1)}})}\\
  &\leq \|w\|_{L^{\infty }(\Omega _{\tau \eta ^{(k_0+1)}})}+\|v-\bar{l}\|_{L^{\infty }(\Omega _{\tau \eta ^{(k_0+1)}})}\\
  &\leq \hat{C}M\eta ^{(k_0+1)(1+\alpha)}.
\end{aligned}
\end{equation*}
Let $l_{k_0+1}=l_{k_0}+\bar{l}$. Recall\cref{e.19} and\cref{e.23}. Then the conclusion holds for $k=k_0+1$. \qed~\\

Now, we can derive the pointwise $C^{1,\alpha}$ regularity for oblique derivative problems in the general form.
\begin{theorem}\label{th5.3}
Let $u$ be a viscosity solution of\cref{e1.1}, $x_0\in \Gamma$ such that $\mathrm{dist}(x_0,\partial\Omega\backslash \Gamma)>1$ and $0<\alpha <\min(\alpha _0,\alpha _1)$. Suppose that $\beta$, $\gamma$ and $g$ are $C^{\alpha}$ at $x_0\in \Gamma$, and $f$ satisfies
\begin{equation*}
\left (\frac{1}{|B_r(x_0)\cap\Omega|}\int_{B_r(x_0)\cap\Omega}|f|^n\right )^{\frac{1}{n}}\leq C_f r^{\alpha - 1}~~\forall r>0.
\end{equation*}
Then $u$ is $C^{1,\alpha}$ at $x_0$, i.e., there exists an affine function $l$ such that
\begin{equation}\label{e11.15}
   \|u-l\|_{L^{\infty }(\bar{\Omega }\cap B_r(x_0))}\leq C_0r^{1+\alpha} ~~~~\forall 0<r<r_1,
\end{equation}
\begin{equation}\label{e41.49}
  |Dl|\leq C_0,
\end{equation}
\begin{equation}\label{e41.50}
  C_0\leq C\left(\|u\|_{L^{\infty }(\Omega )}+C_f+\|g\|_{C^{\alpha }(x_0)}+|F(0)|\right)
\end{equation}
and
\begin{equation}\label{e41.51}
  r_1=\check{C}^{-1},
\end{equation}
where $C$ depends only on $n$, $\lambda$, $\Lambda$, $\delta_0$, $\alpha$ and $\|\gamma \|_{C^{\alpha }(x_0)}$, and $\check{C}$ depends also on $\|\beta \|_{C^{\alpha }(x_0)}$ and the $C^{1}$ modulus of $\Gamma$ at $x_0$.
\end{theorem}
\proof  Rewrite the equation as
  \begin{equation*}
    \left\{
    \begin{aligned}
      &F(D^2u)=f &&\mbox{in}~~  \Omega; \\
      &\beta  \cdot Du= g-\gamma u &&\mbox{on}~~\Gamma .
    \end{aligned}
    \right.
  \end{equation*}
From \Cref{th_ca.1}, $u$ is $C^{\alpha_0}$ at $x_0$. Then, from \Cref{thh4.1}, we obtain that $u$ is $C^{1,\alpha}$ at $x_0$ and\crefrange{e11.15}{e41.51} hold. \qed~\\

Combining with the interior $C^{1,\alpha}$ estimate (see \cite[Theorem 8.3]{C-C}), the boundary local $C^{1,\alpha} $ estimate (\Cref{co4.1}) follows easily (see the proof of \cite[Proposition 2.4]{M-S}).

\section{\texorpdfstring{$C^{2,\alpha}$}{C2,alpha} and higher regularity}\label{S6}
In this section, we prove the $C^{2,\alpha}$ regularity and higher regularity for the oblique derivative problem\cref{e1.1}. We introduce the following two lemmas first for constructing auxiliary functions.
\begin{lemma}\label{le4.1}
Let $F$ be convex and $u$ be a viscosity solution of
\begin{equation*}
    \left\{
    \begin{aligned}
      &F(D^2u)=0 &&\mbox{in}~~ B^+_{1,h_0}; \\
      &\beta \cdot Du = g &&\mbox{on}~~T_1,
    \end{aligned}
    \right.
\end{equation*}
where $\beta $ is a constant vector.

Let $0<\alpha <\alpha _2$. Suppose that $\|g\|_{L^{\infty }(T_r)}\leq C_g r^{1+\alpha}$ for any $0<r<1$. Then there exists a paraboloid $P$ such that
\begin{equation}\label{e6.1}
   \|u-P\|_{L^{\infty }(\bar{B}^{+ }_{r/2,h_0r/2})}\leq C_0r^{2+\alpha}~~\forall 0<r<1,
\end{equation}
\begin{equation}\label{e6.4}
  |DP(0)|+\|D^2P\|\leq C_0,
\end{equation}
and
\begin{equation}\label{e6.5}
  C_0\leq C\left(\|u\|_{L^{\infty }(B^+_{1,h_0})}+ C_g+|F(0)|\right),
\end{equation}
where $C$ depends only on $n$, $\lambda$, $\Lambda$, $\delta_0$, $\alpha$ and $h_0$.
\end{lemma}
\proof  The assumption $F(0)=0$ is made as in \Cref{thh4.1}. Let $M:=\|u\|_{L^{\infty }(B^+_{1,h_0})}+C_g$ and $\Omega_{r}:=B^+_{r,h_0 r}$. We only need to prove the following:

There exist constants $0< \eta < 1$, $\bar{C}$ and $\hat{C}$ depending only on $n$, $\lambda$, $\Lambda$, $\delta_0$, $\alpha$ and $h_0$, and a sequence paraboloids $P_k(x)=\frac{1}{2}x^{T}a_kx+b_kx+c_k$ ($k\geq -1$) such that

\begin{equation}\label{e5.2}
\|u-P_k\|_{L^{\infty }(\Omega _{\eta ^{k}})}\leq \hat{C} M \eta ^{k(2+\alpha )},~~\forall k\geq 0,
\end{equation}

\begin{equation}\label{e5.3}
\eta ^{2k}|a_k-a_{k-1}|+\eta ^{k}|b_k-b_{k-1}|+|c_k-c_{k-1}|\leq \bar{C}\hat{C}M\eta ^{k(2+\alpha)}
\end{equation}
and
\begin{equation}\label{e5.1}
F(a_k)=0, \beta \cdot b_k=0 \mbox{ and } a_k\beta=0.
\end{equation}

We prove it by induction. For $k=0$, by setting $P_0=P_{-1}=0$, the conclusion holds clearly. Suppose that the conclusion holds for $k=k_0$. We need to prove that the conclusion holds for $k=k_0+1$.

Let $r:=\frac{1}{2}\eta ^{k_{0}}$. By \Cref{th2.7}, there exists a unique solution $v$ of
\begin{equation*}
    \left\{
    \begin{aligned}
     & F(D^2v+a_{k_0})=0 &&\mbox{in}~~ B^{+}_{r,h_0r} ; \\
     &\beta  \cdot Dv = 0 &&\mbox{on}~~T_r;\\
     &       v= u-P_{k_0} &&\mbox{on}~~\partial B^{+}_{r,h_0r}\backslash T_r.
    \end{aligned}
    \right.
\end{equation*}
Let $w=u-P_{k_0}-v$. Recall that $\beta\cdot (a_{k_0}x+b_{k_0})=0$ on $T_1$ (see\cref{e5.1}). Then $w$ satisfies
\begin{equation*}
    \left\{
    \begin{aligned}
      &w\in S(\lambda /n,\Lambda ) &&\mbox{in}~~ B^{+}_{r,h_0r}; \\
      &\beta\cdot  Dw = g &&\mbox{on}~~T_r;\\
            &w= 0 &&\mbox{on}~~\partial B^{+}_{r,h_0r}\backslash T_r.
    \end{aligned}
    \right.
\end{equation*}

In the following, we estimate $v$ and $w$ respectively. By the $C^{2,\alpha_2}$ estimate for $v$ (\Cref{th3.10}) and the A-B-P maximum principle (\Cref{th2.4}), there exists a paraboloid
\begin{equation*}
    \bar{P}(x)=\frac{1}{2}x^T\bar{a}x+\bar{b}x+\bar{c}
\end{equation*}
such that
\begin{equation}\label{e5.8}
 \begin{aligned}
 \|v-\bar{P}\|_{L^{\infty }(\Omega _{2\eta r})}\leq &C\frac{(2\eta r)^{2+ \alpha_1 }}{r^{2+ \alpha _{1}}}\|v\|_{L^{\infty }( \Omega _{ r})}\\
 \leq &C\eta ^{2+ \alpha _1}\|u-P_{k_0}\|_{L^{\infty }(\Omega _{ r})}\\
 \leq &C_1\eta ^{\alpha _1-\alpha }\cdot \hat{C}M\eta ^{(k_0+1)(2+\alpha)},
 \end{aligned}
\end{equation}

\begin{equation}\label{e5.4}
 \eta ^{2k_0}|\bar{a}|+ \eta ^{k_0}|\bar{b}|+|\bar{c}|\leq \bar{C}\hat{C}M\eta ^{k_0(2+\alpha)}
\end{equation}
and
\begin{equation}\label{e5.5}
  F(\bar{a}+a_{k_0})=0 \mbox{ and } \beta^0 \cdot \bar{b}=0.
\end{equation}
Furthermore, by the Taylor's formula, for any $(x',0)\in T_r$, we have
\begin{equation*}
  Dv(x',0)=Dv(0)+D^2v(0)(x',0)+o(|x'|).
\end{equation*}
Combining with $\beta \cdot Dv = 0$ on $T_r$, we deduce
\begin{equation}\label{e5.6}
 \bar{a}\beta =0.
\end{equation}

For $w$, by the A-B-P maximum principle, we have
\begin{equation}\label{e5.7}
  \begin{aligned}
\|w\|_{L^{\infty }(\Omega _{r})} \leq Cr\|g\|_{L^{\infty }(T_r)}\leq CC_gr^{2+\alpha}\leq \frac{C_2}{\hat{C}}\hat{C}M\eta ^{k_0(2+\alpha)}.
  \end{aligned}
\end{equation}
Take $\eta $ small enough such that $ C_1\eta ^{\alpha _1- \alpha }<1/2$. Next, take $\hat{C}$ large enough such that
\begin{equation*}
  \frac{C_2}{\hat{C} \eta ^{2+\alpha}}\leq \frac{1}{2}.
\end{equation*}
Therefore, combining\cref{e5.8} and\cref{e5.7}, we have
\begin{equation*}
  \|u-P_{k_0}-\bar{P}\|_{L^{\infty }(\Omega _{\eta ^{(k_0+1)}})}\leq \|w\|_{L^{\infty }}+\|v-\bar{P}\|_{L^{\infty }}\leq \hat{C}M\eta ^{(k_0+1)(2+\alpha)}.
\end{equation*}
Let $P_{k_0+1}=P_{k_0}+\bar{P}$ and recall\crefrange{e5.4}{e5.6}. Then the conclusion holds for $k=k_0+1$. \qed~\\

Based on above $C^{2,\alpha}$ estimate, we deduce the following existence result which will be used to construct auxiliary functions.

\begin{lemma}\label{le4.2}
Let $F$ be convex and $0<\alpha <\min(\alpha_1,\alpha _2)$. Then there exists a unique solution $u\in C^{2,\alpha}(B^{+}_{1,h_0}\cup T_1)\cap C(\bar{B}^{+}_{1,h_0})$ of
\begin{equation}\label{e5.9}
    \left\{
    \begin{aligned}
      &F(D^2u)=0 &&\mbox{in}~~B^{+ }_{1,h_0}; \\
      &(\beta+Ax) \cdot Du =0 &&\mbox{on}~~T_1;\\
      &u=\varphi &&\mbox{on}~~\partial B^{+ }_{1,h_0}\backslash T_1,
    \end{aligned}
    \right.
\end{equation}
where $\beta$ is a constant vector, $\|A\|\leq \delta_0/2$ and $\varphi\in C(\partial B^{+ }_{1,h_0}\backslash T_1)$. Furthermore, we have the estimate
\begin{equation}\label{e6.3}
   \|u\|_{C^{2,\alpha }(\bar{B}^{+ }_{1/2,h_0/2})}\leq C\|\varphi \|_{L^{\infty }(\partial B^{+ }_{1,h_0}\backslash T_1)},
\end{equation}
where $C$ depends only on $n$, $\lambda$, $\Lambda$, $\delta_0$, $\alpha$ and $h_0$.
\end{lemma}

\proof  The existence and uniqueness of $u$ is assured by \Cref{th.existence}. We only need to prove the $C^{2,\alpha}$ regularity. Proving that $u$ is $C^{2,\alpha}$ at $0$ is sufficient, i.e., there exists a paraboloid $P$ such that for any $0<r<1$,
\begin{equation*}
   \|u-P\|_{L^{\infty }(\bar{B}^{+ }_{r/2,h_0r/2})}\leq C\|\varphi \|_{L^{\infty }(\partial B^{+ }_{1,h_0}\backslash T_1)}r^{2+\alpha}.
\end{equation*}

By \Cref{thh4.1}, $u\in C^{1,\alpha}(\bar{B}^{+ }_{1/2,h_0/2})$. It is easy to find a symmetric matrix $B$ such that $B\beta=A^{T}Du(0)$ and $\|B\|\leq C(n,\delta_0)|Du(0)|$. Let $v(x)=u(x)+x^{T}Bx/2$. Then $v$ satisfies
\begin{equation*}
    \left\{
    \begin{aligned}
      &F(D^2v-B)=0 &&\mbox{in}~~ B^+_{1,h_0}; \\
      &\beta \cdot Dv = -Du^{T}Ax+Du(0)^{T}Ax:=g &&\mbox{on}~~T_1,
    \end{aligned}
    \right.
\end{equation*}
Note that $\|g\|_{L^{\infty }(T_r)}\leq C[Du]_{C^{\alpha}} r^{1+\alpha}$ for any $0<r<1$. By \Cref{le4.1}, there  exists a paraboloid $\bar{P}$ such that for any $0<r<1$,
\begin{equation*}
   \|v-\bar{P}\|_{L^{\infty }(\bar{B}^{+ }_{r/2,h_0r/2})}\leq C\left(\|u\|_{L^{\infty }(B^+_{1,h_0})}+ [Du]_{C^{\alpha}(0)}+|F(-B)|\right)r^{2+\alpha}.
\end{equation*}
By the $C^{1,\alpha}$ estimate for $u$,
\begin{equation*}
  [Du]_{C^{\alpha}(0)}+|F(-B)|\leq C \|u\|_{L^{\infty }(B^+_{1,h_0})}.
\end{equation*}
Therefore, let $P=\bar{P}+ x^{T}Bx/2$ and we have
\begin{equation*}
   \|u - P\|_{L^{\infty }(\bar{B}^{+ }_{r/2,h_0r/2})}=\|v-\bar{P}\|_{L^{\infty }}\leq C\|u\|_{L^{\infty }}r^{2+\alpha} \leq C\|\varphi \|_{L^{\infty }(\partial B^{+ }_{1,h_0}\backslash T_1)}r^{2+\alpha}
\end{equation*}
by the A-B-P maximum principle. \qed~\\

Next, we prove the $C^{2,\alpha}$ regularity.

\begin{lemma}\label{th6.1}
Let $F$ be convex and $u$ be a viscosity solution of
  \begin{equation*}
    \left\{
    \begin{aligned}
      &F(D^2u)=f &&\mbox{in}~~ \Omega; \\
      &\beta  \cdot Du = g &&\mbox{on}~~\Gamma,
    \end{aligned}
    \right.
  \end{equation*}
$x_0\in \Gamma$ such that $\mathrm{dist}(x_0,\partial\Omega\backslash \Gamma)>1$ and $0<\alpha <\min(\alpha_1,\alpha _2)$.

Suppose that $\beta$, $g$ and $\Gamma$ are $C^{1,\alpha}$ at $x_0$ and $f$ is $C^{\alpha}$ at $x_0$. Then $u$ is $C^{2,\alpha}$ at $x_0$, i.e., there exists a paraboloid $P$ such that
\begin{equation}\label{e6.6}
   \|u-P\|_{L^{\infty }(\bar{\Omega }\cap B_r(x_0))}\leq C_0r^{2+\alpha} ~~~~\forall 0<r<r_1,
\end{equation}
\begin{equation}\label{e6.7}
  |DP(x_0)|+\|D^2P(x_0)\|\leq C_0,
\end{equation}
\begin{equation}\label{e6.8}
  C_0\leq C\left(\|u\|_{L^{\infty }(\Omega )}+\|f\|_{C^{\alpha } (x_0 )}+\|g\|_{C^{1,\alpha }(x_0)}+|F(0)|\right)
\end{equation}
and
\begin{equation}\label{e6.9}
  r_1=\check{C}^{-1},
\end{equation}
where $C$ depends only on $n$, $\lambda$, $\Lambda$, $\delta_0$ and $\alpha$, and $\check{C}$ depends also on $\|\beta \|_{C^{1,\alpha }(x_0)}$ and $\|\Gamma\|_{C^{1,\alpha}(x_0 )}$.
\end{lemma}
\proof  Similar to the $C^{1,\alpha}$ estimate, we make some normalization first. By choosing a proper coordinate system and scaling, we assume that $x_0$ is the origin,
\begin{equation*}\label{e4.20}
\Gamma  \cap B_1 =\{(x',x_n)\in B_1\big|x_n=\varphi (x')\}
\end{equation*}
and $\varphi $ satisfies that
\begin{equation*}\label{e4.21}
\varphi(0)=0,~~D\varphi(0)=0~~\mbox{and}~~ |\varphi (x')|\leq |x'|^{1+\alpha}.
\end{equation*}
Also, since $\beta$ is defined on $\Gamma  \cap B_1$, we may write $\beta(x')=\beta(x',\varphi(x'))$ and assume that
\begin{equation*}\label{e4.b}
  |\beta(x')-\beta(0)-D\beta(0)x'|\leq |x'|^{1+\alpha} \mbox{ on } \Gamma  \cap B_1
\end{equation*}
and $\|D\beta(0)\|\leq \delta_0/2$. We also assume that $f(0)=0$. Otherwise, we may consider $G(D^2u)=F(D^2u)-f(0)=f-f(0)$. The assumptions that $F(0)=0$, $\beta_n\geq \delta_0 /2$ on $\Gamma  \cap B_1$ are made as in \Cref{thh4.1}.  Similar to $\beta$, we may write $g(x')=g(x',\varphi(x'))$. As in \Cref{thh4.1}, we assume that $g(0)=0$. Furthermore, we assume that $Dg(0)=0$. Otherwise, we may consider
\begin{equation*}
  \tilde{u}(x)=u(x)-\sum_{\beta^0 _i\neq 0,i<n}\frac{g_i(0)x_i^2}{2\beta^0_i}-\sum_{\beta^0_i=0,i<n}\frac{g_i(0)x_ix_n}{\beta^0_n},
\end{equation*}
where $\beta^0:=\beta(0)$. Then $\tilde{u}$ satisfies
 \begin{equation*}
    \left\{
    \begin{aligned}
      &F(D^2\tilde{u}+A)=f &&\mbox{in}~~ \Omega  ; \\
      &\beta \cdot D\tilde{u} = \tilde{g} &&\mbox{on}~~\Gamma  ,
    \end{aligned}
    \right.
  \end{equation*}
where $A$ is a constant matrix and $\tilde{g}$ satisfies $D\tilde{g}(0)=0$. Then the desired estimates for $u$ follow easily from that of $\tilde{u}$.

Let $M:=\|u\|_{L^{\infty }(\Omega  )}+\|f\|_{C^{\alpha }( 0 )}+\|g\|_{C^{1, \alpha }(0)}$ and $\Omega _{r}:=\Omega \cap B_{r,h_0 r}$. We only need to prove the following:

There exist constants $0< \tau <1$, $0< \eta < 1$, $\bar{C}$ and $\hat{C}$ depending only on $n$, $\lambda$, $\Lambda$, $\delta_0$ and $\alpha$, and a sequence paraboloids $P_k(x)=\frac{1}{2}x^{T}a_kx+b_kx+c_k$ ($k\geq -1$) such that

\begin{equation}\label{e4.24}
\|u-P_k\|_{L^{\infty }(\Omega _{\tau \eta^{k}})}\leq \hat{C} M \eta ^{k(2+\alpha )},~~\forall k\geq 0,
\end{equation}

\begin{equation}\label{e4.25}
\eta ^{2k}|a_k-a_{k-1}|+\eta ^{k}|b_k-b_{k-1}|+|c_k-c_{k-1}|\leq \bar{C}\hat{C}M\eta ^{k(2+\alpha )}
\end{equation}
and
\begin{equation}\label{e4.36}
F(a_k)=0, \beta^0 \cdot b_k=0 \mbox{ and } \sum_{j=1}^{n}\left((a_k)_{ij}\beta ^{0}_j+D_i\beta _{j}(0)(b_k)_j\right)=0,~~i<n.
\end{equation}

We prove it by induction. For $k=0$, by setting $P_0=P_{-1}=0$, the conclusion holds clearly. Suppose that the conclusion holds for $k=k_0$. We need to prove that the conclusion holds for $k=k_0+1$.

Let $r:=\frac{1}{2}\tau \eta ^{k_{0}}$, $\tilde{B}^{+}_{r,h_0r}:=B^{+}_{r,h_0r}- r^{1+\alpha} e_n $, $\tilde{T}_r:=T_r- r^{1+\alpha} e_n$ and $\tilde{\Omega }_{r}:=\Omega \cap \tilde{B}^{+}_{r,h_0r}$. Assume that
\begin{equation}\label{e4.37}
  4r^{\alpha}\leq h_0.
\end{equation}
Then $\Omega _{r/2}\subset \tilde{\Omega }_{r}$.

Note that $u-P_{k_0}$ satisfies
\begin{equation*}
  F(D^2(u-P_{k_0})+a_{k_0})=f~~\mbox{in}~~\Omega _{2r}.
\end{equation*}
By the H\"{o}lder estimate (\Cref{th5.1}),
\begin{equation}\label{e4.28}
\begin{aligned}
  \|u-P_{k_0}\|&_{L^{\infty }(\Omega _{r})}+r^{\alpha _0}[u-P_{k_0}]_{C^{\alpha _0}(\Omega _{r})}\\
  \leq &C\left (\|u-P_{k_0}\|_{L^{\infty }(\Omega _{2r})}+r^{2}\|f\|_{L^{\infty }(\Omega _{2r})}+r\|g\|_{L^{\infty }(\Gamma \cap B_{2r})}\right ).
\end{aligned}
\end{equation}
Extend $u-P_{k_0}$ from $\Omega_r$ to the whole $R^n$ such that
\begin{equation}\label{e4.23}
  \|u-P_{k_0}\|_{L^{\infty }(R^n)}+r^{\alpha_0}[u-P_{k_0}]_{C^{\alpha_0}(R^n )}\leq C \left ( \|u-P_{k_0}\|_{L^{\infty }(\Omega_{r})}+r^{\alpha_0}[u-P_{k_0}]_{C^{\alpha_0}(\bar{\Omega} _{r})}\right ).
\end{equation}

Let $\alpha<\alpha_3<\min(\alpha_1,\alpha_2)$. By \Cref{le4.2}, there exists a unique solution $v\in C^{2,\alpha_3}$ of
\begin{equation*}
    \left\{
    \begin{aligned}
     & F(D^2v+a_{k_0})=0 &&\mbox{in}~~ \tilde{B}^{+}_{r,h_0r} ; \\
     & (\beta^0+D\beta (0)x') \cdot Dv = 0 &&\mbox{on}~~\tilde{T}_r;\\
     &       v= u-P_{k_0} &&\mbox{on}~~\partial \tilde{B}^{+}_{r,h_0r}\backslash\tilde{T}_r.
    \end{aligned}
    \right.
\end{equation*}
Let $w=u-P_{k_0}-v$. Then $w$ satisfies
\begin{equation*}
    \left\{
    \begin{aligned}
      &w\in S(\lambda /n,\Lambda , f) &&\mbox{in}~~ \Omega \cap \tilde{B}^{+}_{r,h_0r}; \\
      &\beta\cdot  Dw = g-\beta\cdot (a_{k_0}x+b_{k_0}+Dv) &&\mbox{on}~~\Gamma  \cap \tilde{B}^{+}_{r,h_0r};\\
            &w= 0 &&\mbox{on}~~\partial \tilde{B}^{+}_{r,h_0r}\cap \bar{\Omega}.
    \end{aligned}
    \right.
\end{equation*}

In the following arguments, we estimate $v$ and $w$ respectively. By the boundary $C^{2,\alpha_3}$ estimates for $v$ (\Cref{le4.2}) and the A-B-P maximum principle (\Cref{th2.4}), there exists a paraboloid
  \begin{equation*}
    \tilde{P}(x)=\frac{1}{2}\left(x+ r^{1+\alpha}e_n\right)^T\bar{a}\left(x+r^{1+\alpha}e_n\right)+\bar{b}\left(x+r^{1+\alpha}e_n\right)+\bar{c}  \end{equation*}
such that
\begin{equation*}
\begin{aligned}
  \|v-\tilde{P}\|_{L^{\infty }(\Omega _{2\eta r})}\leq& C\frac{(2\eta r)^{2+ \alpha_3 }}{r^{2+ \alpha _{3}}}\|v\|_{L^{\infty }( \tilde{B}^{+}_{r,h_0r})}\leq C\eta ^{2+ \alpha _3}\|u-P_{k_0}\|_{L^{\infty }(\Omega _{ r})} \\
  \leq& C\eta ^{\alpha _3-\alpha }\cdot \hat{C}M\eta ^{(k_0+1)(2+\alpha)},
\end{aligned}
\end{equation*}
\begin{equation}\label{e4.26}
 \eta ^{2k_0}|\bar{a}|+ \eta ^{k_0}|\bar{b}|+|\bar{c}|\leq \bar{C}\cdot \hat{C}M\eta ^{k_0(2+\alpha)}
\end{equation}
and
\begin{equation}\label{e4.38}
  F(\bar{a}+a_{k_0})=0 \mbox{ and } \beta^0 \cdot \bar{b}=0.
\end{equation}
Furthermore, by the Taylor's formula, for any $(x',-r^{1+\alpha})\in \tilde{T}_r$, we have
\begin{equation*}
  Dv(x',-r^{1+\alpha})=Dv(0,-r^{1+\alpha})+D^2v(0,-r^{1+\alpha})(x',0)+o(|x'|).
\end{equation*}
Combining with $(\beta^0+D\beta (0)x') \cdot Dv = 0$ on $\tilde{T}_r$, we deduce
\begin{equation}\label{e4.39}
 \sum_{j=1}^{n}\left(\bar{a}_{ij}\beta ^{0}_j+D_i\beta _{j}(0)\bar{b}_j\right)=0,~~i<n.
\end{equation}

Let $\bar{P}(x)=\frac{1}{2}x^T\bar{a}x+\bar{b}x+\bar{c}$. Then
\begin{equation}\label{e4.27}
\begin{aligned}
   \|v-\bar{P}\|_{L^{\infty }(\Omega _{2\eta r})}\leq & \|v-\tilde{P}\|_{L^{\infty }(\Omega _{2\eta r})}+ \|\bar{P}-\tilde{P}\|_{L^{\infty }(\Omega _{2\eta r})}\\
   \leq & C\eta ^{\alpha _3-\alpha }\cdot \hat{C}M\eta ^{(k_0+1)(2+\alpha)}+|Kr^{1+\alpha}e_n^T\bar{a}x|\\
   &+|K^2r^{2+2\alpha}\bar{a}_{nn}|+|Kr^{1+\alpha} \bar{b}_n|\\
  \leq & \left (C_1\eta ^{\alpha _3-\alpha }+\frac{C_2 \tau^{1+\alpha}}{\eta ^{2+\alpha }} \right )\cdot \hat{C}M\eta ^{(k_0+1)(2+\alpha)}.
\end{aligned}
\end{equation}

Next, we estimate the term $w$. Let $\tilde{B}_{\mu }:=\{x\in \tilde{B}^{+ }_{r,h_0r}\big|\mathrm{dist}(x,\partial \tilde{B}^{+ }_{r,h_0r}\backslash\tilde{T}_r)\geq \mu r\}$, $\Omega _{r,\mu }:=\Omega \cap \tilde{B}_{\mu }$, $\Gamma_1:=\partial \tilde{B}_{\mu }\cap \bar{\Omega} $ and $\Gamma_2:=\tilde{B}_{\mu }\cap \Gamma $.

By the global H\"{o}lder estimate for $v$ (\Cref{th2.7}) and recalling\cref{e4.28} and\cref{e4.23}, there exists $0<\alpha_4\leq \alpha_0/2$ such that
\begin{equation}\label{e4.29}
\begin{aligned}
    \|v\|_{L^{\infty }(\tilde{\Omega} _{r})}+r^{\alpha _4}[v]_{C^{\alpha _4}(\tilde{\Omega} _{r})}\leq & C\left (\|u-P_{k_0}\|_{L^{\infty }(\Omega _{r})}+r^{\alpha _0}[u-P_{k_0}]_{C^{\alpha _0}(\Omega _{r})}\right )\\
    \leq &C\left(\|u-P_{k_0}\|_{L^{\infty }(\Omega _{2r})}+r^{2}\|f\|_{L^{\infty }(\Omega _{2r})}+r\|g\|_{L^{\infty }(\Gamma \cap B_{2r})}\right). \\
\end{aligned}
\end{equation}
For any $ x_0\in \Gamma_1$, there exists $ \bar{x}\in \partial \tilde{B}^{+ }_{r,h_0r}\backslash \tilde{T}_r$ such that $|x_0-\bar{x}|=\mu r$. Then by recalling\cref{e4.28}
\begin{equation}\label{e4.30}
\begin{aligned}
|w(x_0)|=&|u(x_0)-P_{k_0}(x_0)-v(x_0)|\\
     =&|u(x_0)-P_{k_0}(x_0)-v(x_0)-u(\bar{x})+P_{k_0}(\bar{x})+v(\bar{x})| \\
    \leq & |\left (u(x_0)-P_{k_0}(x_0)\right )-\left (u(\bar{x})-P_{k_0}(\bar{x})\right )|+|v(x_0)-v(\bar{x})|\\
    \leq &\frac{C(\mu r)^{\alpha _4}}{r^{\alpha _4}}\left(\|u-P_{k_0}\|_{L^{\infty }(\Omega _{2r })}+r^{2}\|f\|_{L^{\infty }(\Omega _{2r})}+r\|g\|_{L^{\infty }(\Gamma \cap B_{2r})}\right).
\end{aligned}
\end{equation}
By the A-B-P maximum principle, combining with\cref{e1.28}, we have
\begin{equation}\label{e4.31}
\begin{aligned}
\|w\|&_{L^{\infty }(\Omega _{r,\mu })} \\
\leq &\|w\|_{L^{\infty }(\Gamma_1)}+Cr\|g-\beta\cdot  (a_{k_0}x+b_{k_0}+Dv)\|_{L^{\infty }(\Gamma_2)}+Cr^{2}\|f\|_{L^{\infty }(\Omega _{2r})}\\
\leq &\frac{C(\mu r)^{\alpha _4}}{r^{\alpha _4}}\left(\|u-P_{k_0}\|_{L^{\infty }(\Omega _{2r})}+Cr^{2}\|f\|_{L^{\infty }(\Omega _{2r})}
+r\|g\|_{L^{\infty }(\Gamma \cap B_{2r})}\right)\\
&+Cr\|g\|_{L^{\infty }(\Gamma \cap B_{2r})}+Cr\|\beta\cdot(a_{k_0}x+b_{k_0}+Dv)\|_{L^{\infty }(\Gamma_2)}+Cr^{2}\|f\|_{L^{\infty }(\Omega _{2r})}\\
\leq& C_3\mu ^{\alpha _4}\hat{C}M\eta ^{k_0(2+\alpha)}+C_4r^{2+\alpha}\|g\|_{C^{1,\alpha }(0)}+C_5r^{2+\alpha}\|f\|_{C^{\alpha }(0)}\\
&+Cr\|\beta\cdot(a_{k_0}x+b_{k_0})\|_{L^{\infty }(\Gamma_2)}+Cr\|\beta\cdot  Dv\|_{L^{\infty }(\Gamma_2)}\\
\leq & \left (C_3\mu ^{\alpha _4}+\frac{C_4}{\hat{C}}+\frac{C_5}{\hat{C}}\right )\hat{C}M\eta ^{k_0(2+\alpha)}\\
&+Cr\|\beta\cdot(a_{k_0}x+b_{k_0})\|_{L^{\infty }(\Gamma_2)}+Cr\|\beta\cdot  Dv\|_{L^{\infty }(\Gamma_2)}.
  \end{aligned}
\end{equation}

In the following, we estimate $\|\beta\cdot(a_{k_0}x+b_{k_0})\|_{L^{\infty }(\Gamma_2)}$ and $\|\beta\cdot  Dv\|_{L^{\infty }(\Gamma_2)}$ respectively. For the first term, recall\cref{e4.25} and\cref{e4.36} and then we obtain
\begin{equation}\label{e4.32}
  \begin{aligned}
  \|\beta&\cdot(a_{k_0}x+b_{k_0})\|_{L^{\infty }(\Gamma_2)}=\|(\beta-\beta^{0}-D\beta(0)x')\cdot  (a_{k_0}x+b_{k_0})\|_{L^{\infty }(\Gamma_2)}\\
  &+\|(\beta^{0}+D\beta(0)x')\cdot  (a_{k_0}x+b_{k_0})\|_{L^{\infty }(\Gamma_2)}\\
\leq& r^{1+\alpha}\|a_{k_0}x+b_{k_0}\|_{L^{\infty }(\Gamma_2)}+|D\beta(0)||a_{k_0}|r^2+\|\beta^0\cdot a_{k_0}x+b_{k_0}\cdot D\beta(0)x'\|_{L^{\infty }(\Gamma_2)}\\
\leq& \frac{C\hat{C}M}{1-\eta }r^{1+\alpha }+\frac{C\hat{C}M}{1-\eta }r^2+\|\sum_{j=1}^{n}(a_{k_0})_{nj}\beta ^{0}_jx_n\|_{L^{\infty }(\Gamma_2)}\\
\leq& C_6\hat{C}Mr^{1+\alpha }
\end{aligned}
\end{equation}

We assume that
\begin{equation}\label{e4.33}
  4 \tau^{\alpha} \leq \mu.
\end{equation}
Then, $\forall x_0 \in \Gamma_2$, $\mathrm{dist} (x_0,\tilde{T}_r )< \frac{1}{2}\mathrm{dist} (x_0,\partial \tilde{B}^{+}_{r,h_0r}\backslash \tilde{T}_r )$. Let $x^{\ast }\in \tilde{T}_r$ such that $|x_0-x^{\ast}|=\mathrm{dist} (x_0,\tilde{T}_r )$. By the $C^{2,\alpha _3}$ estimate for $v$ in $B_{\mu r}(x^{\ast})\cap \tilde{B}^{+}_{r,h_0r}$ and noting that $(\beta^{0}+D\beta(0)x^{'}_0)\cdot Dv(x^{\ast})=0$, we have
\begin{equation}\label{e4.34}
  \begin{aligned}
   |\beta&(x_0) \cdot Dv(x_0)|\\
   \leq & |\beta(x_0) -\beta^{0}-D\beta(0)x^{'}_0||Dv(x_0)|+|\beta^{0}+D\beta(0)x^{'}_0||Dv(x_0)-Dv(x^{\ast})|\\
   \leq &  r^{1+\alpha }\frac{C}{\mu r}\|v\|_{L^{\infty }(B_{\mu r}(x^{\ast})\cap \tilde{B}^{+}_{r,h_0r})}+C\frac{2 r^{1+\alpha }}{(\mu r)^{2}}\|v\|_{L^{\infty }(B_{\mu r}(x^{\ast})\cap \tilde{B}^{+}_{r,h_0r})}\\
   \leq & \frac{C_7 r^{\alpha}}{\mu}\cdot \hat{C}M\eta ^{k_0(2+\alpha)}+\frac{C_8r^{\alpha-1}}{\mu^{2}}\cdot \hat{C}M\eta ^{k_0(2+\alpha)}\\
  \end{aligned}
\end{equation}
Combining\cref{e4.31},\cref{e4.32} and\cref{e4.34}, we have
\begin{equation}\label{e4.35}
  \|w\|_{L^{\infty }(\Omega _{r,\mu })}\leq \left(C_3\mu ^{\alpha _4}+\frac{C_4+C_5}{\hat{C}}+C_6\tau ^{2+ \alpha }+\frac{C_7\tau ^{1+\alpha }}{\mu }+\frac{C_8\tau ^{\alpha } }{\mu^2 }\right)\hat{C}M\eta ^{k_0(2+\alpha)}.
\end{equation}

Take $\eta $ small enough such that $ C_1\eta ^{\alpha _3- \alpha }<1/4$. Let $\mu =\tau^{\alpha /4}$ and take $\tau$ small enough such that
\begin{equation*}
\begin{aligned}
   &\frac{C_2\tau^{1+\alpha }}{\eta ^{2+ \alpha }}<\frac{1}{4} ,~\frac{C_3\tau ^{\alpha \alpha _4/4}}{\eta ^{2+ \alpha }}<\frac{1}{12},~\frac{C_6\tau ^{2+\alpha }}{\eta ^{2+\alpha}}<\frac{1}{12}, \frac{C_7 \tau ^{1+3\alpha/4 }}{\eta ^{2+\alpha}}<\frac{1}{12},\frac{C_8\tau ^{\alpha/2 } }{\eta ^{2+ \alpha }}<\frac{1}{12}
\end{aligned}
\end{equation*}
and\cref{e4.37} and\cref{e4.33} hold. Finally, take $\hat{C}$ large enough such that
\begin{equation*}
  \frac{C_4+C_5}{\hat{C} \eta ^{2+\alpha}}\leq \frac{1}{12}.
\end{equation*}
Therefore, combining\cref{e4.27} and\cref{e4.35}, we have
\begin{equation*}
\begin{aligned}
  \|u-P_{k_0}-\bar{P}\|_{L^{\infty }(\Omega _{\tau \eta ^{(k_0+1)}})}= & \|u-P_{k_0}-v+v-\bar{P}\|_{L^{\infty }(\Omega _{\tau \eta ^{(k_0+1)}})}\\
  \leq& \|w\|_{L^{\infty }(\Omega _{\tau \eta ^{(k_0+1)}})}+\|v-\bar{P}\|_{L^{\infty }(\Omega _{\tau \eta ^{(k_0+1)}})}\\
  \leq &\hat{C}M\eta ^{(k_0+1)(2+\alpha)}.
\end{aligned}
\end{equation*}
Let $P_{k_0+1}=P_{k_0}+\bar{P}$; Recall\cref{e4.26},\cref{e4.38} and\cref{e4.39}. Then the conclusion holds for $k=k_0+1$. \qed~\\

Similar to the pointwise $C^{1,\alpha}$ estimate, we have the following:
\begin{theorem}\label{th5.333}
Let $F$ be convex, $u$ be a viscosity solution of\cref{e1.1} and $x_0\in \Gamma$ such that $\mathrm{dist}(x_0,\partial\Omega\backslash \Gamma)>1$ and $0<\alpha <\min(\alpha_0,\alpha_1,\alpha _2)$. Suppose that $\beta$, $g$ and $\Gamma$ are $C^{1,\alpha}$ at $x_0$ and $f$ is $C^{\alpha}$ at $x_0$.

Then $u$ is $C^{2,\alpha}$ at $x_0$, i.e., there exists a paraboloid $P$ such that
\begin{equation*}\label{e6.666}
   \|u-P\|_{L^{\infty }(\bar{\Omega }\cap B_r(x_0))}\leq C_0r^{2+\alpha} ~~~~\forall 0<r<r_1,
\end{equation*}
\begin{equation*}\label{e6.777}
  |DP(x_0)|+\|D^2P(x_0)\|\leq C_0,
\end{equation*}
\begin{equation*}\label{e6.888}
  C_0\leq C\left(\|u\|_{L^{\infty }(\Omega )}+\|f\|_{C^{\alpha } (x_0 )}+\|g\|_{C^{1,\alpha }(x_0)}+|F(0)|\right)
\end{equation*}
and
\begin{equation*}\label{e6.999}
  r_1=\check{C}^{-1},
\end{equation*}
where $C$ depends only on $n$, $\lambda$, $\Lambda$, $\delta_0$, $\alpha$ and $\|\gamma \|_{C^{1,\alpha }(x_0)}$, and $\check{C}$ depends also on $\|\beta \|_{C^{1,\alpha }(x_0)}$ and $\|\Gamma\|_{C^{1,\alpha}(x_0 )}$.
\end{theorem}

Combining with the interior $C^{2,\alpha}$ estimate (see \cite[Theorem 8.1]{C-C}), the boundary local $C^{2,\alpha} $ estimate (\Cref{th4.6}) follows (see the proof of \cite[Proposition 2.4]{M-S}).

Since we have obtained the $C^{2,\alpha}$ regularity, the higher regularity for the oblique derivative problem can be deduced standardly.
\begin{theorem}\label{th.hig}
Let $F$ be convex, $u$ be a viscosity solution of\cref{e1.1} and $0<\alpha <\min(\alpha_0,\alpha_1,\alpha _2)$. Suppose that $F\in C^{k,\alpha}(S^n)$, $\Gamma\in C^{k+2,\alpha} $, $\beta,\gamma,g\in C^{k+1,\alpha }(\bar{\Gamma}  )$ and $f\in C^{k,\alpha }(\bar{\Omega}  )$.

Then for any $\Omega' \subset \subset \Omega\cup \Gamma$, $u\in C^{k+2,\alpha}(\bar{\Omega '})$ and
\begin{equation}\label{eq.hig}
   \|u\|_{C^{k+2,\alpha }(\bar{\Omega '})}\leq C,
\end{equation}
where $C$ depends on $n$, $\lambda$, $\Lambda$, $\delta_0$, $\alpha $, $\|\beta \|_{C^{k+1,\alpha }(\bar{\Gamma}  )}$, $\|\gamma  \|_{C^{k+1,\alpha }(\bar{\Gamma}  )}$, $\|f\|_{C^{k,\alpha }(\bar{\Omega })}$, $\|g\|_{C^{k+1,\alpha }(\bar{\Gamma}  )}$, $\|u\|_{L^{\infty }(\Omega )}$, $F$, $\Omega '$ and $\Omega$.

In particular, if $F$, $\Gamma , f,\beta, \gamma, g \in C^{\infty }$, then $u\in C^{\infty }(\bar{\Omega' })$.
\end{theorem}
\proof  $u\in C^{k+2,\alpha}(\Omega )$ is well known and we only need to prove the boundary $ C^{k+2,\alpha}$ estimate. We prove the theorem by induction. Let $k=1$. For any $x_0\in \Gamma  $, there exists a proper coordinate system such that $x_0$ is the origin,
\begin{equation*}\label{e4.45}
\Gamma  \cap B_2 =\{(x',x_n)\in B_2\big|x_n=\varphi (x')\},
\end{equation*}
where $\varphi\in C^{k+2,\alpha} (T_2) $ satisfies
\begin{equation*}\label{e4.46}
\varphi(0)=0,~~D\varphi(0)=0.
\end{equation*}
Introduce the transformation $y=\psi (x)$ where $\psi $ is defined as follows: $y_i=x_i$ for $i<n$ and $y_n=x_n-\varphi (x')$. Define $\tilde{u}(y)=u(x)$,
\begin{equation*}
 A_{ij}=\frac{\partial y_i}{\partial x_j} \mbox{ and } B_{ij}=u_l \cdot \frac{\partial x_l}{\partial y_m} \cdot \frac{\partial y_m}{\partial x_i\partial x_j}.
\end{equation*}
It is easy to check that $\tilde{u}$ is a viscosity solution of
\begin{equation}\label{e4.48}
\left\{
\begin{aligned}
      &G(D^2\tilde{u},y)=\tilde{f} &&\mbox{in}~~ B^+_{1,h_0}; \\
      &\tilde{\beta} \cdot D\tilde{u} +\tilde{\gamma} \tilde{u}= \tilde{g} &&\mbox{on}~~T_1,
\end{aligned}
\right.
\end{equation}
where $G(M,y):=F(A^{T}MA+B)$ for any $M\in S^n$, $\tilde{f}(y):=f(x)$, $\tilde{\beta _j}(y):=\beta _i (x)\partial y_j/\partial x_i$, $\tilde{\gamma }(y):=\gamma (x) $ and $\tilde{g}(y):=g(x)$.

By \Cref{th4.6}, $u\in C^{2,\alpha}(B_1\cap \bar{\Omega})$. Hence, $A_{ij},B_{ij}\in C^{1,\alpha}$. Combining with $F\in C^{1,\alpha}$, we have that $G\in C^{1,\alpha}(S^n\times \bar{B}^+_{1,h_0})$. Note that $\tilde{\beta _n}(y)=\beta _i\cdot  \partial y_n/\partial x_i=\beta_n(x)-\beta'(x)\cdot D\varphi (x')$. Since $D\varphi(0)=0$, by a proper scaling, we may assume that
\begin{equation*}
  \tilde{\beta _n}\geq \frac{\delta_0 }{2} \mbox{ on } T_1.
\end{equation*}
Then $\tilde{u}\in C^{3,\alpha}(B^+_{1,h_0})$ can be obtained by the classical interior estimates (see \cite [Proposition 9.1]{C-C} ). Differentiate the equations with respect to $x_m$ ($m\leq n-1$) and we have
\begin{equation}\label{e4.44}
    \left\{
    \begin{aligned}
      &G_{ij}(D^2\tilde{u}(y),y)(\tilde{u}_m)_{ij}=-G_m(D^2\tilde{u}(y),y)+\tilde{f}_m(y) &&\mbox{in}~~ B^+_{1,h_0}; \\
      &\tilde{\beta} \cdot D\tilde{u}_m +\tilde{\gamma} \tilde{u}_m= g-\tilde{\beta} _{m}\cdot D\tilde{u}-\tilde{\gamma}_m\tilde{u} &&\mbox{on}~~T_1.
    \end{aligned}
    \right.
\end{equation}
By the boundary estimates for linear elliptic equations (see \cite [Theorem 4.40] {Li7}), we have that $\tilde{u}_m\in C^{2,\alpha}(\bar{B} ^+_{3/4,3h_0/4})$ and
\begin{equation*}
     \|\tilde{u}_m\|_{C^{2,\alpha }(\bar{B} ^+_{3/4,3h_0/4})}\leq C.
\end{equation*}
Thus, $\tilde{u}\in C^{3,\alpha}(T_{3/4})$. From the boundary estimates for the Dirichlet problems, we have that $\tilde{u}\in C^{3,\alpha}(\bar{B} ^+_{1/2,h_0/2})$ and
\begin{equation*}
     \|\tilde{u}\|_{C^{3,\alpha }(\bar{B} ^+_{1/2,h_0/2})}\leq C.
\end{equation*}
Hence, $u\in C^{3,\alpha }(B_{1/2}\cap \bar{\Omega })$. Since $x_0\in \Gamma  $ is arbitrary, by the standard covering argument, we have that $u\in C^{3,\alpha }(\bar{\Omega' })$ and the estimate\cref{eq.hig} holds.

Assume that the theorem holds for $k=k_0$. We prove that the theorem hods for $k=k_0+1$. Since $u\in C^{k_0+2,\alpha }(\bar{\Omega '})$, $B_{ij}\in C^{k_0+1,\alpha }$ where $B_{ij}$ is defined as above.  From \cite [Proposition 9.1]{C-C} we know that $u\in C^{k_0+3,\alpha}(\Omega \cap B_1)$ and hence $\tilde{u}\in C^{k_0+3,\alpha}(B^+_{1,h_0})$. Differentiate\cref{e4.48} $k_0+1$ times with respect to the horizontal directions. Then we deduce equations similar to\cref{e4.44}. From the regularity for linear oblique derivative problems, we obtain that the $k_0+1$ order horizontal derivatives of $\tilde{u}$ lie in $C^{2,\alpha}(\bar{B}^+_{3/4,3h_0/4})$. Hence, $\tilde{u}\in C^{k_0+3,\alpha}(T_{3/4})$. From the regularity for Dirichlet problems, we obtain that $\tilde{u}\in C^{k_0+3,\alpha}(B^+_{1/2,h_0/2})$ and hence $u\in C^{k_0+3,\alpha }(B_r\cap \bar{\Omega })$ for some $r>0$. Therefore, by a scaling and covering argument, $u\in C^{k_0+3,\alpha }(\bar{\Omega' })$ and the estimate\cref{eq.hig} holds. \qed~\\

\begin{acknowledgements}
The authors would like to thank Professor Xinan Ma and Professor Yu Yuan for their helps. This research was supported by NSFC 11671316 and NSFC 11701454.
\end{acknowledgements}

~\\
\textbf{Compliance with ethical standards}~\\
\textbf{Conflict of interest}~\\
Both authors declare that they have no potential conflict of interest.

\end{document}